\documentclass[11pt,twoside]{amsart}
\usepackage{amsmath}
\usepackage{amssymb, verbatim}
\usepackage{graphicx}
\usepackage{comment}
\usepackage{bbm}
\usepackage{mathtools}
\usepackage[utf8]{inputenc}
\usepackage[english]{babel}
\usepackage[usenames, dvipsnames]{color}
\usepackage[margin= 1in]{geometry}
\usepackage{enumitem}
\usepackage{booktabs}

\title[]{Intermediate disorder regime for half-space directed polymers}
\author{Xuan Wu}
\address{Department of Mathematics, University of Chicago, Eckhart Hall
5734 S University Ave
Chicago IL, 60637}
\email{xuanw@uchicago.edu}

\numberwithin{equation}{section}
\theoremstyle{plain}
\newtheorem{thm}{Theorem}[section]

\newtheorem{defn}[thm]{Definition}
\newtheorem{lem}[thm]{Lemma}

\theoremstyle{definition}

\newtheorem{rk}[thm]{Remark}

\newcommand{\del}{\partial}
\newcommand\norm[1]{\left\lVert#1\right\rVert}

\newcommand{\bt}{{\bf{t}}}
\newcommand{\bi}{{\bf{i}}}
\newcommand{\bx}{{\bf{x}}}
\newcommand{\vi}{\mathbf{i}}
\newcommand{\vx}{\mathbf{x}}

\newcommand{\vt}{\mathbf{t}}

\DeclareMathOperator{\Var}{Var}
\DeclareMathOperator{\sym}{Sym}

\DeclarePairedDelimiter\floor{\lfloor}{\rfloor}
 
\newcommand{\N}{\mathbb{N}}
\newcommand{\R}{\mathbb{R}}
\newcommand{\Z}{\mathbb{Z}}
\newcommand{\Pro}{\mathbb{P}}
\newcommand{\E}{\mathbb{E}}

\allowdisplaybreaks

\renewcommand{\epsilon}{\varepsilon}
\newcommand{\ov}[1]{\overline{#1}}

\begin{document}
\maketitle
\begin{abstract}
We consider the convergence of point-to-point partition functions for the half-space directed polymer model in dimension 1+1 in the intermediate disorder regime as introduced for the full space model by Alberts, Khanin and Quastel in \cite{AKQ}. By scaling the inverse temperature as $\beta n^{-1/4}$, the point-to-point partition function converges to the chaos series for the solution to stochastic heat equation with Robin boundary condition and delta initial data. Furthermore, the convergence result is then applied to the exact-solvable log-gamma directed polymer model in a half-space.
\end{abstract}

\section{Introduction}
The directed polymers were introduced in the statistical physics literature by Huse and Henley \cite{HH} in 1985 and received first rigorous mathematical treatment in 1988 by Imbrie and Spencer \cite{IS}. The monograph \cite{Com} is a great resource for the foundational work in this area. Over the last thirty
years, the directed polymers played an important role as a playground of many fascinating problems in the probability world. 

Among those different directions opened up by directed polymers, in dimension 1+1, its connection to the KPZ universality class \cite{Cor} has attracted extensive attention. The polymer measure in dimension 1+1 is a random probability measure on paths in a random environment, which favors higher weighted paths. It is constructed through up / right paths on $\mathbb{Z}^2$ with path measure re-weighted by an i.i.d. random environment presented at each lattice points. The KPZ universality conjecture concerns the large scale asymptotic behavior of the polymer free energy and there are two characteristic scalings , the 1:2:3 KPZ scaling and the weak noise scaling, known as the strong KPZ universality conjecture and the weak KPZ universality conjecture respectively.

In the direction of the strong KPZ universality conjecture for directed polymers, the first rigorous verification of the 1/3 fluctuation of polymer free energy was proven for a special case \cite{Sep}, where the integrable log-gamma polymers were introduced. Among directed polymers, the log-gamma directed polymer model was special in the same way as the last passage percolation models with exponential or geometric weights are special among corner growth models. Namely, both demonstrate integrable structures
and permit explicit computations. \cite{COSZ} computed the Laplace transform of the point-to-point partition function. \cite{BCR} transformed that formula into a Fredholm determinant and performed asymptotic analysis, with motivation from Macdonald process formulas in \cite{BC}.

Under the weak noise scaling, the convergence of polymer free energy in dimension 1+1 to KPZ equation has been established in the remarkable work by Alberts, Khanin and Quastel in \cite{AKQ}, which is known to have proved the weak KPZ universality conjecture for directed polymers.

It is natural to ask the same question for the half-space polymers. The half-space directed polymers are
constructed through up / right paths constrained to stay in the half-quadrant with path measure re-weighted by two random environments($X$ present only at the boundary and $\omega$ in the bulk). Compared to the full space case, the extra boundary environment $X$ penalizes or rewards the path measure every time the walker visits the origin in an i.i.d. manner. The {\bf{main Theorem~\ref{main_thm}}} of this paper builds the connection between half-space directed polymers and half-space stochastic heat equation(SHE) with Robin boundary condition/KPZ equation with Neumann boundary condition. 

Aside from the general half-space polymer model, recently there has also been considerable attention focused on the exact-solvable log-gamma polymers, see the recent work in \cite{BBC} and \cite{OSZ}. But presently no rigorous asymptotics have been proved. This motivates to apply the convergence results for general half-space polymer model to the log-gamma case, see Section~\ref{application-section}. Our result was further used in \cite{Par1} to obtain an equality-in-distribution for SHE on the half space with different boundary conditions. 

More generally, half-space KPZ universality is also studied by other half-space models approached from the perspective of scaling to KPZ equation and also from the perspective of exact solvability. On half-space asymmetric simple exclusion process (ASEP), \cite{CS} showed that the height function
converges to Hopf-Cole solution of KPZ equation with Neumann boundary
condition(Robin boundary condition for SHE). With stronger estimates developed, \cite{Par2} extended their results to negative values of the boundary condition. In the exact
solvability direction, \cite{BBCW} studied half-line ASEP as a scaling limit of a
stochastic six-vertex model in a half-quadrant and found exact formulas for
half-space KPZ/SHE with $\mu = -1/2$, see \eqref{SHE}. See also in \cite{GH} for the study of KPZ
equation with Neumann boundary conditions in the context of the theory of
regularity structures.

\subsection*{Outline} In Section~\ref{formulation} we give a precise formulation of our main result Theorem~\ref{main_thm} and heuristics of the proof are provided in Section~\ref{heuristics}. The techniques we borrow from U-statistics are stated in Section~\ref{U-stats}. Our main technical estimates are provided in Section~\ref{heat_kernel_estimates} with proofs postponed to the appendix. We leave the proof of our main theorem to Section~\ref{proof-section}. In the last Section~\ref{application-section}, we discuss the half-space log-gamma polymer model and apply our main theorem to get an analogous convergence result for the point-to-point partition function.

\subsection*{Acknowledgement}
The author is very grateful to Ivan Corwin for his incredible guidance and many encouraging conversations, and also extends thanks to Guillaume Barraquand and Promit Ghosal for their helpful discussions related to this work. The author is in particular very grateful to an anonymous referee for pointing out many typos/errors and for providing numerous suggestions. The author was partially supported by  NSF grant of Ivan Corwin's DMS-1664650 as well as the NSF grant DMS-1441467 for PCMI, at which this work started.

\section{Definitions of the model and main results}\label{formulation} 
The aim of this paper is to study the SHE limits of half-space directed polymers in a random environment. We start with definitions of the half-space polymers.

\subsection{Half-space Polymers}\label{def-half-space-polymer}
Consider an $n$-step simple symmetric random walk on non-negative integers $\mathbb{N}_{ 0}$ with a totally reflecting barrier at the origin. The law of this walk is equal to that of the absolute value of a standard symmetric random walk on $\mathbb{Z}$. Denote the reflecting random walk probability measure by $\Pro_{R}$ on paths starting from origin at time $0$ and we also denote $\Pro^{m,x}_{R}$ as the probability measure on paths starting at $x\geq 0$ at time $m\geq 0$.  This measure $\Pro_{R}^{m,x}$ will serve as our background probability measure throughout this paper and we omit the superscript when there is no ambiguity about the starting point and time. For a path $S$, let $S_i$ denote its location at time $i$ and define transition probability for a random walk starting at $x$ at time $m$ and arriving at $y\geq 0$ at time $n\geq m$ by
\[
p(m,n,x,y) \coloneqq \sum_{S: S_n = y} \Pro^{m,x}_R(S).
\]

Such path measures will be affected by two environments and we start with the boundary environment. Let $X=\{X_i\}$ be a sequence of i.i.d. non-negative random variables and we refer to $X$ as the boundary random environment. Define the random transition kernel as
\begin{align}\label{transition_reflect}
p_X(m,n,x,y) \coloneqq \sum_{S:S_n=y} \left(\prod_{m\leq i < n :S_i= 0}X_i\right)\cdot \Pro^{m,x}_{R}(S).
\end{align}

Denote $\N$ as the set of positive natural numbers while $\N_0$ also includes zero and denote $[m,n]_{\Z}$ as the integers inside $[m,n]$. Given a path $S:[m,n]\to \mathbb{N}_0$, define the corresponding random measure $\Pro_X$ as
\[
\Pro_X(S) \coloneqq \left(\prod_{ m\leq i<n:S_i= 0}X_i\right)\Pro_{R}(S).
\]
$\Pro_X$ is a measure-valued random variable with randomness inherited from $X$. Note that in general $\Pro_X$ is not a probability measure due to the ``punishing" or ``rewarding" effects caused by the random environment $X$ when paths visit the origin.

When the boundary random environment is deterministic such that $X_i \equiv \gamma \geq 0$, $\gamma$ is denoted as the reflection rate for the barrier at origin. It follows that the barrier is absorbing if $0\leq \gamma < 1$, totally reflecting if $\gamma = 1$, and rewarding if $\gamma >1$. Now the transition kernel $p_{\gamma}(m,n,x,y)$ also becomes deterministic. Explicitly,
\begin{align}\label{def:pgamma}
p_{\gamma}(m,n,x,y)\coloneqq \sum_{j=0}^{n-m}\gamma^j\mathbb{P}^{m,x}_R(N_{m,n}=j,S_n=y).
\end{align}
Here $N_{m,n}$ is the total visits to the origin as
\begin{align}\label{def:localtime}
N_{m,n}(S)\coloneqq \#\{i\in [m,n-1]_\mathbb{Z}\ |\ S_i=0\}.
\end{align}

Let $\omega(i,x)$ for $(i,x) \in \mathbb{N}_{ 0} \times \mathbb{N}_{ 0}$ be an i.i.d. collection of random variables and we refer to $\omega \coloneqq \{\omega(i,x)\}$ as the bulk random environment. The energy of an $n$-step nearest neighbor walk $S$ in the environment $\omega$ is defined as:
\[
H^{\omega}_n(S) \coloneqq \sum_{i=0}^{n-1} \omega(i,S_i).
\]
Define the polymer probability measure with randomness inherited from both the bulk random environment $\omega$ and the boundary random environment $X$ as:
\begin{align*}
\Pro^{\omega,X}_{n,\beta}(S) &\coloneqq \frac{1}{Z^{\omega,X}(n;\beta)}e^{\beta H^{\omega}_n(S)}\cdot \Pro_X(S)\\
&=\frac{1}{Z^{\omega,X}(n;\beta)}e^{\beta H^{\omega}_n(S)}\cdot \left(\prod_{ 0 \leq i < n:S_i = 0}X_i\right)\cdot\Pro_{R}(S).
\end{align*}
Here $\beta$ is a parameter, called inverse temperature. The normalization term $Z^{\omega,X}(n;\beta)$ is a point-to-line partition function, defined as:
\[
Z ^{\omega,X}(n;\beta) \coloneqq \E_R\left[e^{\beta H^{\omega}_n(S)}\left(\prod_{0 \leq i < n :S_i = 0}X_i\right)\right],
\]
where the expectation is taken with respect to the reflecting random walk measure $\Pro_R$ and preserves randomness from $\omega$ and $X$.

The main {\bf goal} of this paper is to study the limiting behavior of the following point-to-point partition function: 
\begin{equation}\label{p2p-part-function}
Z^{\omega,X} \left(n,x; \beta\right) \coloneqq \E_R \left[ e^{\beta H^{\omega}_n(S)}\left(\prod_{ 0 \leq i < n:S_i = 0}X_i\right)\cdot \mathbbm{1} \{ S_n = x \} \right],
\end{equation}
where $\mathbbm{1}$ is the indicator function. Note that $\{S_n=x\}$ is non-empty only if $n$ and $x$ have the same parity, which we denote as $n\leftrightarrow x$. Generally, for $n\in\mathbb{N}$ and $x\in\mathbb{R}$, denote $[x]_n$ as the largest integer among which are smaller than $x$ and enjoys the same parity as $n$, i.e.
\begin{align}\label{def:parity}
[x]_n\coloneqq\max\{m\in\mathbb{Z}\ |\ m\leq x,\ m\leftrightarrow n\}.
\end{align} 

\subsection{Stochastic Heat Equation with Robin boundary condition}
In this section we introduce the SHE with Robin boundary condition, which arises as a weak scaling limit of the half-space directed polymers. We also provide the expression of the chaos series for its solution,  a series of multiple stochastic integrals over a Robin heat kernel with respect to a space-time white noise. 

\subsubsection{1-D heat equation with Robin boundary condition}
\begin{defn}
We say $\rho_{\mu}(t,x,y)$ is the fundamental solution to 1-D heat equation on $\mathbb{R}_{\geq 0}$ with Robin boundary condition and initial data $\delta(y-x)$ if
\begin{align}\label{Heat-eq}
\del_t \rho_{\mu}(t,x,y) &=\frac{1}{2} \del_{xx} \rho_{\mu}(t,x,y)\\
\del_x \rho_{\mu}|_{x=0} &= \mu \cdot\rho_{\mu}|_{x=0},\nonumber
\end{align}
and if for any function $\varphi(x)$,
\[
v(t,x) = \int_0^{\infty} \rho_{\mu}(t,x,y)\varphi(y) dy
\]
solves heat equation with initial condition
\[
v(0,x) = \varphi(x).
\]
\end{defn}
There are a few equivalent forms of the Robin heat kernel. We will make use of the following form which can be found in \cite[Lemma 4.4]{CS}.
\begin{align} \label{Robin-kernel}
\rho_{\mu}(t,x,y)= &( {2\pi t})^{-1/2} \left( e^{ - (y-x)^2/(2t)} - e^{ - (y+x)^2/(2t)}  \right)\\
&+ 2( 2\pi t^3)^{-1/2}  \int_0^{\infty}(y+x+s) e^{-\mu s- (y+x+s)^2/(2t) } ds. \nonumber
\end{align}

\subsubsection{Stochastic Heat equation with Robin boundary condition}
Consider the stochastic heat equation with multiplicative noise
\begin{equation}\label{SHE}
\del_t z_{\beta} = \frac{1}{2}\del_{xx} z_{\beta} + \beta z_{\beta}\cdot \xi
\end{equation}
with delta initial data and Robin boundary condition:
\begin{align*}
z_{\beta}(0,\cdot) &= \delta(0)\\
\del_x z_{\beta}(\cdot,x)|_{x=0} &= \mu\cdot z_{\beta}(\cdot,0).
\end{align*}
Here $\xi(t,x)$ is a white noise on $\R_{\geq 0}\times\R_{\geq 0}$ with covariance structure
\[
\E[ \xi(t,x) \xi(s,y) ] = \delta(t-s) \delta(x-y).
\]
For details about white noise and full space SHE, we refer to \cite[Section 3]{AKQ}. Further discussions can be found in \cite{Jan}.

With the help of the Robin heat kernel, the mild solution is given by
\begin{align}\label{mild_solution}
z_{\beta}(t,x) = \sum_{k=0}^{\infty}\int_{\Delta_k (t)}\int_{\mathbb{R}_{\geq 0}^k} \rho_{\mu}(t-t_k,x_k,x)\cdot\beta^k \prod_{i=1}^k \rho_{\mu}(t_i -t_{i-1}, x_{i-1}, x_i)d\xi^{\otimes k}(\bt,\bx),
\end{align}
where $\Delta_k(t) = \{0=t_0 < t_1 < \cdots < t_k \leq t \}$ and $x_0=0$. 

To simply the notation, we define the k-fold operator as follows. Let $k\in\mathbb{N}_0$ and $g(t_1,t_2,x_1,x_2)$ be a function defined on $0\leq t_1<t_2$ and $(x_1,x_2)\in\mathbb{R}^2$.  $\mathbf{F}_k[g](t,x;\bt,\bx):(\mathbb{R}_{>0} \times \mathbb{R})\times \Delta_k(t)\times \mathbb{R}^k\to \mathbb{R}$ is defined as
\begin{equation}\label{def:Fk}
\mathbf{F}_k[g](t,x;\bt,\bx)\coloneqq g(t_{k},t,x_k,x )\prod_{j=1}^k g(t_{j-1},t_{j},x_{j-1},x_{j}).
\end{equation}
Here  the convention $t_0=x_0=0$ has been used. Let 
\begin{align}
\rho_{\mu,k }(t,x;\bt,\bx)=\mathbf{F}_k[\rho_\mu](t,x;\bt,\bx),
\end{align}
with the understanding that $\rho_\mu(s,t,\cdot,\cdot) \coloneqq \rho_\mu(t-s,\cdot,\cdot)$. Then
\begin{align*} 
z_{\beta}(t,x) = \sum_{k=0}^{\infty}\int_{\Delta_k (t)}\int_{\mathbb{R}_{\geq 0}^k}\beta^k  \rho_{\mu,k }(t,x;\bt,\bx)  d\xi^{\otimes k}(\bt,\bx).
\end{align*}

Our main result below shows that by diffusively scaling the random walks, under intermediate disorder scaling($\beta n^{-1/4}$) and critical scaling near the boundary, the point-to-point partition function converges to $z_{\sqrt{2}\beta}(t,x)$, solution to SHE. The convergence takes place in the topology of supremum norm on bounded continuous functions, denoted as $\xrightarrow{(d)}$. Denote $\lambda(\beta) = \log \E[e^{\beta \omega}]$, our main theorem is as follows.
\begin{thm}\label{main_thm}
Fix $\mu\in\mathbb{R}$. Let $\omega$ be i.i.d. random environment with mean zero and variance one which satisfies $\lambda(\beta_0) < \infty$ for some $\beta_0>0$. For $n\in\mathbb{N}$, let $\gamma=1-\mu/\sqrt{n}$. Assume that $X$ satisfies $\E[X]=\gamma$ and that $\E\left[|X-\E[X]\right|^3] \leq Kn^{-\epsilon} $ for some $\epsilon\in (0,1]$ and $K>0$.  Then
\begin{equation*}
2^{-1}n^{1/2}e^{-\floor*{nt} \lambda(\beta n^{-1/4})}  Z^{\omega,X}\bigl(\floor*{nt}, [x\sqrt {n}]_{\floor*{nt}}; 
\beta n^{-{1}/{4}} \bigr) \mathop{\longrightarrow}^{(d)} z_{\sqrt{2}\beta}(t,x).
\end{equation*}
Here $[x\sqrt{n}]_{\floor*{nt}}$ is the largest integer which is smaller than $x\sqrt{n}$ and has the same parity as $\floor*{nt}$. See \eqref{def:parity}.
\end{thm}

\begin{rk}
Here we require the third moment assumption in order to prove tightness and we do not believe this is the optimal case.
\end{rk}

\section{Heuristics and ideas of proof}\label{heuristics}
In this section we attempt to explain why $\beta n^{-1/4}$ and $\gamma = 1- {\mu}/{\sqrt{n}}$ are natural scalings. We also provide heuristics behind the proof of the Main Theorem~\ref{main_thm} and comment on the main technical ingredients. First let us summarize the setup of half-space polymers in the following Table~\ref{halfspace}. Note that in the left picture, random walk trajectories are pictured as paths in a half-quadrant while the partition functions are defined with respect to random walks on non-negative integers. The equivalence between these two formulations is clear and in this way the figure better illustrate the idea. For simplicity of notations, we omit the floor function when it does not cause ambiguity, e.g. $\floor*{nt},  [x\sqrt{n}]_{\floor*{nt}}$.
\begin{table}[h]
     \begin{center}
     \begin{tabular}{c p{9cm}}
     \hline
      Half-quadrant polymers with $\omega$ and $X$ & \quad Definition of partition functions   \\ 
    \cmidrule(r){1-2}
     \raisebox{-4cm}{\includegraphics[width= 7cm, height=4cm]{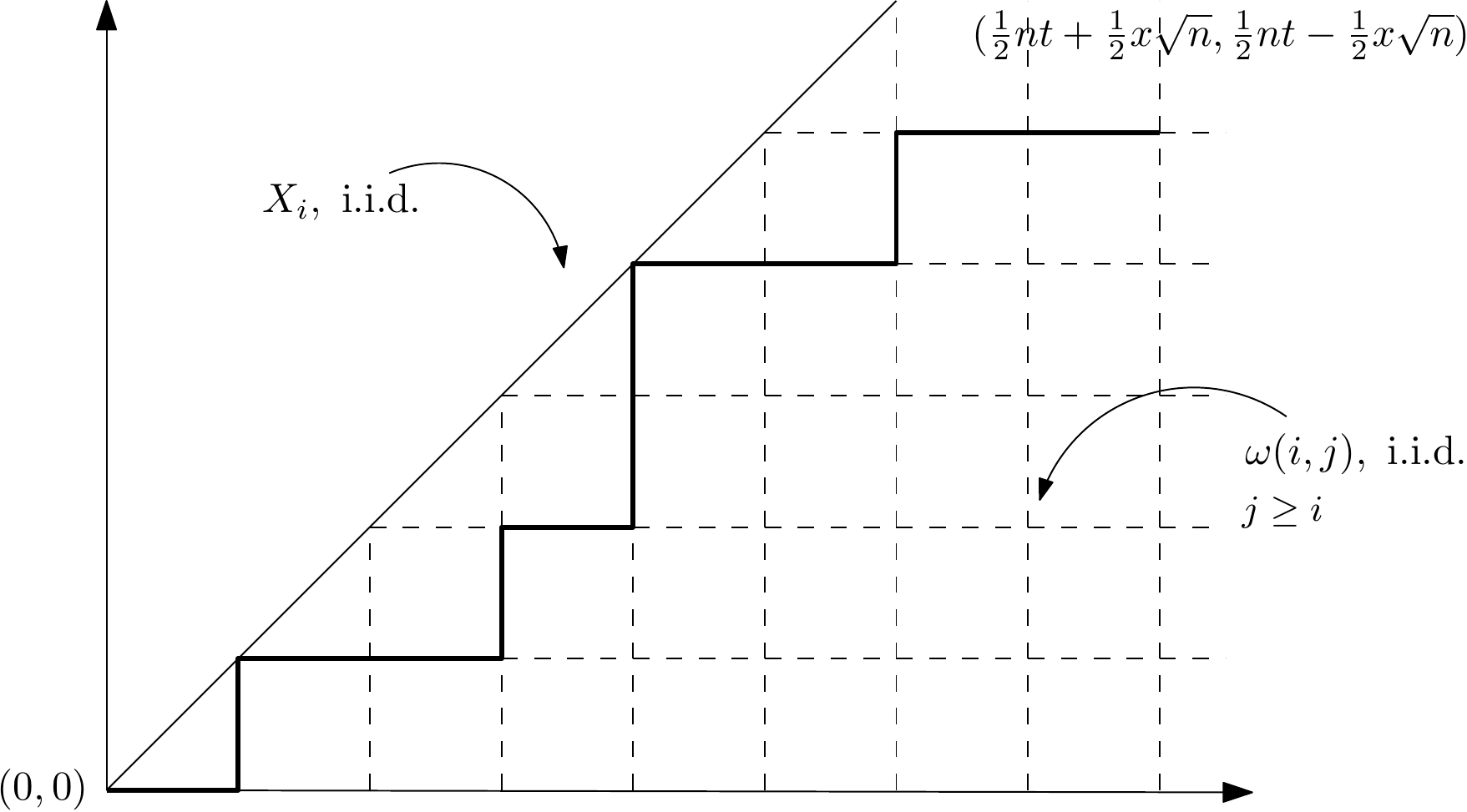}}
      & 
      Let $\beta >0$, the scaled point-to-point partition function is defined as:
      
      \medskip
      
     $\displaystyle{Z^{\omega,X} \left(nt, x\sqrt{n}; \beta n^{-1/4}\right)} :=$
     
     \medskip
     
     $\displaystyle{ \E_R \left[ e^{\beta n^{-1/4} H^{\omega}_{nt}(S)}\left(\prod_{ \substack{0\leq i<nt\\S_i = 0}}X_i\right)\cdot \mathbbm{1} \{ S_{nt} = x\sqrt{n} \} \right]}.$
      
      \\      
      \\
      \hline
      \end{tabular}
      \caption{Summary of the half-space polymer model. The expectation is taken with respect to probability measure $\Pro_R$ and preserves randomness from $\omega$ and $X$. $H^{\omega}_n(S) := \sum_{i=0}^{n-1} \omega(i,S_i)$ is the energy of an $n$-step nearest neighbor walk $S$}.
      \label{halfspace}
      \end{center}
      \end{table}

The tuning at boundary, $\gamma = 1- {\mu}/{\sqrt{n}}$, is clear. When the background random walk is scaled diffusively, the total number of visits to the boundary of this random walk is of scale $\sqrt{n}$. In the average sense, in order to see a non-trivial limit of $\displaystyle{\prod_{0\leq i<n:S_i = 0}X_i}$, we must have $\gamma-1=O\left(\frac{1}{\sqrt{n}}\right)$.

The strategy for proving Theorem~\ref{main_thm} is to first prove the convergence for a modified partition function $\mathfrak{Z}^{\omega}$. $\mathfrak{Z}^{\omega}$ takes the form of a discrete chaos series with random walk transition probability kernel. The techniques of U-statistics in Section~\ref{U-stats} provide criteria for convergence of discrete chaos series to continuous ones. Furthermore we rewrite the unmodified partition function $Z^{\omega}$ in the same form as $\mathfrak{Z}^{\omega_n}$ with a perturbed environment $ {\omega}_n$, still of mean zero but with variance only asymptotically one. In addition, the same strategy will be applied in the log-gamma polymer model, where we will need to deal with the issue that the random environment will only be i.i.d. on the diagonal and the bulk respectively.

Denote $D_k^n$ as a discrete integer simplex:
\begin{align}\label{D_k^n}
D_k^n \coloneqq \{ \bi = (i_1 \cdots, i_k) \in \N_{ 0}^k : 0\leq i_1 < \cdots < i_k < n  \}.
\end{align}
We define a $k$-fold transition kernel $p_{X,k}(n,y;\vi, \mathbf{y})$ for $(n,y,\vi, \mathbf{y}) \in(\mathbb{N}\times \mathbb{N}_0)\times D_k^n \times \N_{0}^k$ of a half-space random walk with a barrier at origin that arrives at $y$ in $n$ steps.
\begin{align}\label{k-fold-prob} 
p_{X,k }(n,y;\vi, \mathbf{y}) & \coloneqq p_X(i_k,n, y_k, y)\prod_{j=1}^k p_X(i_{j-1},i_j ,y_{j-1}, y_j).
\end{align}
Here the convention $i_0=x_0=0$ has been used. The modified point-to-point partition is defined as:
\begin{align}\label{modified-p2p}
\displaystyle{\mathfrak{Z}^{\omega,X} \left(n, y; \beta\right)}&\coloneqq \E_R \left[ \prod_{i=0}^{n-1} \bigl( 1 + \beta
 \omega(i, S_i) \bigr)\cdot\left(\prod_{0 \leq i < n : S_i = 0}X_i\right) \cdot \mathbbm{1} \{ S_{n}= y \} \right]
\end{align}
Expanding the above product in the expectation and by a direct computation, $\mathfrak{Z}^{\omega,X}(n,y; \beta  )$ could be written as a discrete sum of weighted chaos,
\begin{align}\label{modified-p2p-part-function}
\mathfrak{Z} ^{\omega,X}(n,y; \beta )
&= p_X(0,n,0,y) + \sum_{k=1}^{n}
\beta^k  \sum_{\vi\in D_k^n} \sum_{\mathbf{y}\in\N^k_0} p_{X,k }(n,y;\vi, \mathbf{y})
\omega(\bi,\mathbf{y}),
\end{align}
where $\omega(\vi,\mathbf{y}) \coloneqq \displaystyle\prod_{j=1}^k \omega(i_j, y_{j})$.

Heuristically we may see why SHE \eqref{mild_solution} arise in the limit. Under the diffusive scaling and boundary tuning ($\gamma = 1- {\mu}/{\sqrt{n}}$), random walk transition probabilities converge to the Robin heat kernel. Moreover, the random environment $\omega$ approximates White noise by scaling $\beta$ to zero in a critical manner (i.e. $\beta n^{1/4}$). To made this rigorous, we need local limit theorem and $L^2$ bounds on the $k$-fold random transition kernels $p_{X,k}(n,y;\vi, \mathbf{y})$. These are the main technical inputs of these paper and are provided in Section~\ref{heat_kernel_estimates}. 

To see that $\beta n^{-1/4}$ is the critical scaling, it is illustrative to check that the $k=1$ term in the summation above has order $O(\sqrt{n})$. For simplicity, assume $X_i \equiv 1$ and consider the point-to line case, i.e. do not fix the endpoint. Now it suffices to show that $n^{-1/4}\displaystyle{\sum_i\sum_x\omega(i, x)\mathbb{P}(S_i = x)}$ stays bounded as a random variable (with randomness inherited from $\omega$). This could be easily seen from taking the second moment. In detail, we see that
\begin{align*}
&\quad\E_{\omega}\left(n^{-1/4}\sum_i\sum_x\omega(i, x)\mathbb{P}(S_i = x)\right)^2 \\
&= n^{-1/2}\E_{\omega}\left[\sum_i\sum_x\omega(i, x)\mathbb{P}(S_i = x)\sum_j\sum_y\omega(j, y)\mathbb{P}(\tilde{S}_j = y)\right]\\
&=n^{-1/2}\sum_{i,j}\sum_{x,y}\E_{\omega}\left[\omega(i, x)\omega(j, y)\right]\mathbb{P}(S_i = x)\mathbb{P}(\tilde{S}_j = y)\\
&=n^{-1/2}\sum_{i=j}\sum_{x=y}\mathbb{P}(S_i = x)\mathbb{P}(\tilde{S}_j = y)\\
&= O(1).
\end{align*}
Here $S,\tilde{S}$ are two independent random walk paths. The third equality follows from taking expectation with respect to $\omega$ by Fubini theorem. Only the the intersection points of $S,\tilde{S}$ will contribute to the sum as $\omega$ is i.i.d. of mean zero and variance one. From general theory of 1-D random walks, we know that $S$ and $\tilde{S}$ intersect $O(\sqrt{n})$ times on average and this explains the scaling $\beta n^{-1/4}$.

\section{U-statistics}\label{U-stats}
The techniques of U-statistics are convenient for obtaining convergence of partition functions $\mathfrak{Z}^{\omega}$, which take the form of discrete chaos. As the results about U-statistics are already presented in \cite[Section 4]{AKQ}, we choose to state the results and refer the proofs to their counterparts in \cite{AKQ}. See \cite{CSZ} for a more general treatment of discrete chaos expansion with more general random environment.

We start with introducing the definition of U-statistics and then quote a technical lemma (Lemma~\ref{lemma_U-stat}). In application to log-gamma polymer models, we need to allow a slightly more general setting. See Lemma \ref{sum_sk_tilde_omega}.

Recall that $n\leftrightarrow x$ denotes $n$ and $x$ have the same parity.  More generally, $\bi \leftrightarrow \mathbf{y} $ means that all corresponding entries share the same parity. Let $\mathcal{R}_k^n$ be the collection of rectangles, defined as:
\[
\mathcal{R}_k^n\coloneqq \biggl\{ \left[ n^{-1} \vi ,
n^{-1}(\vi+ \mathbf{1}) \right) \times\left[ n^{-1/2} {\mathbf{y}} , n^{-1/2} (\mathbf{y}+\mathbf
{2})  \right): \vi\in
D_k^n, \mathbf{y}\in\mathbb{N}_0^k, \vi\leftrightarrow\mathbf{y}
\biggr\}.
\]
Here $D_k^n$ is integer simplex defined in \eqref{D_k^n} and $\mathbf{1}$ is the $k$-dimensional vector $(1,1,\cdots, 1)$. Also
\[
\left[ n^{-1} \vi ,
n^{-1}(\vi+ \mathbf{1}) \right) \coloneqq \left[ n^{-1}  i_1 ,
n^{-1}(i_1+ 1) \right)\times\cdots\times \left[ n^{-1}  i_k ,
n^{-1}(i_k+ 1) \right),
\]
and similarly,
\[
\left[ n^{-1/2} {\mathbf{y}} , n^{-1/2} (\mathbf{y}+\mathbf
{2})  \right) \coloneqq \left[ n^{-1/2} {y_1} , n^{-1/2} (y_1+2)  \right)\times\cdots\times \left[ n^{-1/2} {y_k} , n^{-1/2} (y_k+2)\right).
\]

For a $L^2$ function $g$ on $[0,1]^k \times\R_{\geq 0}^k$, take $n \geq 1$, the corresponding U-statistics $\mathcal{S}_k^n(g)$ of $g$ could be viewed as a weighted average of a discretization of $g$ through the random environment $\omega$. We now discretize $L^2([0,1]^k \times \R_{\geq 0}^k)$ functions by replacing their values with their integral mean values on rectangles in $\mathcal{R}_k^n$. Consider a function $g\in L^2([0,1]^k \times \R_{\geq 0}^k)$, define $\ov{g}_n$ by specifying the values of $\ov{g}_n$ on every $R \in \mathcal{R}_k^n$, more specifically we define
\[
\ov{g}_n|_{R} \coloneqq \frac{1}{|R|}\int_R g.
\]
where $|R| = 2^kn^{-3k/2}$. Note that $\ov{g}$ is constant on every single $R$ and for each $n,k$ fixed, each pair $(\vi,\mathbf{y})\in D_k^n\times \mathbb{N}_{0}^k$ ($\bi \leftrightarrow \mathbf{y}$) corresponds to a unique $R \in \mathcal{R}_k^n$.

For the convenience of applying U-statistics results we consider sums over unordered sets
\[
E_k^n\coloneqq  \{ \vi\in[1,n]_\mathbb{Z}^k :
\vi_j \neq\vi_l \textrm{ for } j \neq l \}.
\]
Recall that $\omega(\vi,\mathbf{y}) = \displaystyle\prod_{j=1}^k \omega(i_j, y_{j})$.

\begin{defn}
The corresponding U-statistics of $g \in L^2([0,1]^k \times\R_{\geq 0}^k)$ is defined as
\begin{equation}\label{Def-U-stats}
\mathcal{S}_k^n(g;\omega) \coloneqq 2^{k/2}\sum_{\bi\in E^n_k}\sum_{\mathbf{y} \in \mathbb{N}_{0}^k} \ov{g}_n\left(n^{-1} {\bi} ,n^{-1/2} {\mathbf{y}} \right) \omega(\bi,\mathbf{y})\cdot\mathbbm{1}\{\bi \leftrightarrow \mathbf{y}\}.
\end{equation}
\end{defn}

The following lemma, proved as \cite[Lemma 4.1]{AKQ}, bounds the second moment of $\mathcal{S}_k^n(g;\omega)$ from above.
\begin{lem}\label{lemma-U-bounds}
Let $\mathcal{S}_k^n(g;\omega)$ be a U-statistics as in \eqref{Def-U-stats}. For any linear combinations of functions $g_1, \cdots, g_m \in L^2([0,1]^k \times \mathbb{R}_{\geq 0}^k)$ through $\alpha_1, \cdots,\alpha_k \in \mathbb{R}$, we have
\[
\sum_{l=1}^m \alpha_l \mathcal{S}_k^n(g_l;\omega) = \mathcal{S}_k^n\left( \sum_{l=1} ^ m \alpha_l g_l ;\omega\right).
\]

Moreover, if random environment variables satisfy moment conditions $\E[\omega(i,x)] = 0$ and $\Var [\omega(i,x)] =\sigma^2$, then 
\[
\E \bigl[ \mathcal{S}_k^n(g)^2 \bigr] \leq \sigma^{2k} n^{3k/2} 
\norm{g}^2_{L^2([0,1]^k\times \mathbb{R}_{\geq 0}^k)}.
\]
\end{lem}

Note that the U-statistics is invariant under permutation for $(\mathbf{t},\bx)$ and we denote 
\[
\sym g(\mathbf{t},\bx) = \frac{1}{k!}\sum_{\pi\in \sigma_k} g(\pi \mathbf{t},\pi \bx),
\]
where $\sigma_k$ is the symmetric group of degree $k$.

For $G = (g_0, g_1, g_2, \ldots)\in \bigoplus_{k \geq0} L^2([0,1]^k
\times\R_{\geq 0}^k)$, define its chaos series $I(G)$ as follows,
\begin{align*}
I(G)&\coloneqq \sum_{k=0}^{\infty} \int_{[0,1]^k} \int_{\R^k_{\geq 0}} \sym g_k(\vt,
\vx) \xi^{\otimes k}( d\vt\, d \vx)\\
&= \sum_{k=0}^{\infty} \int_{[0,1]^k} \int_{\R^k_{\geq 0}} g_k(\vt,
\vx) \xi^{\otimes k}( d\vt\, d \vx).
\end{align*}

 The following lemma, proved as \cite[Theorem 4.5]{AKQ}, shows that under mild conditions, the U-statistics converges in distribution to the continuum chaos series. 
\begin{lem}\label{lemma_U-stat} 
Let $\omega_n(i,x),\ (i,x)\in\mathbb{N}_0\times \mathbb{N}_{0}$ be a sequence of i.i.d. random environments that satisfy
\begin{align}\label{omega}
 \E[\omega_n]=0,\ \textup{and}\ \lim_{n\to\infty}\E[\omega^2_n]=1.
\end{align}  
Let $G = (g_0, g_1, g_2, \ldots)\in \bigoplus_{k \geq0} L^2([0,1]^k
\times\R_{\geq 0}^k)$ with 
\begin{align}\label{GG}
\lim_{N\to\infty}\limsup_{n\to\infty} \sum_{k=N}^{\infty} \E[\omega^2_n]^k \norm{g_k}_{L^2} =0.
\end{align}
Then as $n \to\infty$,
\[
\mathcal{S}^n(G)\coloneqq \sum_{k=0}^{\infty}
n^{-{3k}/{4}} \mathcal{S}_k^n(g_k;\omega_n) \xrightarrow{{\hspace*{0.4cm}(d)}\hspace*{0.4cm}} I(G).
\]
Moreover, suppose $G_1, \ldots, G_m \in\bigoplus_{k \geq0} L^2([0,1]^k
\times\R_{\geq 0}^k)$ all satisfy \eqref{GG}. Then as $n \to\infty$, we have the joint convergence:
\[
\bigl( \mathcal{S}^n(G_1), \ldots ,\mathcal{S}^n(G_m)
\bigr) \xrightarrow{{\hspace*{0.4cm}(d)}\hspace*{0.4cm}}\bigl(I(G_1), \ldots, I(G_m)
\bigr).
\]
\end{lem}
For the application to log-gamma polymer model, we need the following lemma for a perturbed random environment $\tilde{\omega}$.

\begin{lem}\label{sum_sk_tilde_omega}
Let $\tilde{\omega}_n(i,x),\ (i,x)\in\mathbb{N}^2_0 $ be a sequence of random environments. Assume that for fixed $n$, $\tilde{\omega}_n(i,x),(i,x)\in\mathbb{N}_0 \times \mathbb{N} $ are i.i.d. random variables and that $\tilde{\omega}_n(i,0), i \in \mathbb{N}_0$ are also i.i.d. random variables. Furthermore, assume that $\E[\tilde{\omega}_n(i,0)]=\E[\tilde{\omega}_n(i,1)]=0$, $\lim_{n\to\infty} \E[\tilde{\omega}^2_n(i,1)]=1$ and that
\begin{equation*}
\sigma^2\coloneqq \sup_{n\in\mathbb{N},x\in\mathbb{N}_0}\E[\tilde{\omega}^2_n(i,x)]<\infty.  
\end{equation*}
Then, replacing \eqref{GG} with
\begin{align}\label{tGG}
\lim_{N\to\infty}\limsup_{n\to\infty} \sum_{k=N}^{\infty} \sigma^{2k} \norm{g_k}_{L^2} =0,
\end{align} 
the convergence results in Lemma~\ref{lemma_U-stat} still hold with $\omega_n$ replaced by $\tilde{\omega}_n$.  
\end{lem}

\begin{proof}
The proof follows as a trivial reasoning in \cite[Theorem 4.5]{AKQ}.
\end{proof}

\section{Estimates on discrete transition kernel}\label{heat_kernel_estimates}
We record in this section estimates that will be needed in proving Theorem \ref{main_thm}. Their proofs are postponed to the appendix. Recall that $p_\gamma$ defined in \eqref{def:pgamma} is the deterministic transition kernel as $X\equiv \gamma$. Lemma \ref{transition_p_unform_upperbd} concerns pointwise upper bounds for $p_\gamma$. In particular, it shows that $p_\gamma$ enjoys exponential decay. Lemma \ref{pointwise_convergence_tran_barrier} proves the local limit theorem for $p_\gamma$. Lemma \ref{var_transition} bounds the variance of the random transition kernel $p_X$ in terms of $p_\gamma$. Combined with Lemma \ref{transition_p_unform_upperbd}, it implies the variance of $p_{X}$ also decays exponentially. 

\begin{lem}\label{transition_p_unform_upperbd}
For any $\mu\in\mathbb{R}$ and $\tau\geq 1$, there exist  a constant  $ B_0(\mu,\tau)$ and a universal constant $C_0$ such that the following statement holds. For $n\in\mathbb{N}$, $m\in [1,\tau n]_\mathbb{Z}$ and $(x, y)\in \mathbb{N}^2_0$, let $\gamma=1-\mu/\sqrt{n}$. Then
\begin{align*}
p_\gamma (0,m, x ,y )\leq 2 B_0(\mu,\tau)m^{-1/2}e^{-(x-y)^2/(C_0m)}. 
\end{align*}
\end{lem}
\begin{lem}\label{pointwise_convergence_tran_barrier}
For any $\mu\in\mathbb{R}$, $0<\epsilon< 1$ and $M,\tau\geq 1$, there exists $\textup{Err}(n;\mu,\epsilon,M,\tau)$ such that the following statement holds. Assume that $n\in\mathbb{N}$, $t\in [\epsilon,\tau ]$ and $(x,y)\in [0,M]^2$ with $nt\in\mathbb{N}$, $\sqrt{n}x,\sqrt{n}y\in\mathbb{Z}$ and $nt\leftrightarrow \sqrt{n}(y-x)$. Let $\gamma = 1 - {\mu}/{\sqrt{n}}$. Then
\begin{align*}
\left|2^{-1}n^{1/2} p_\gamma (0,nt ,\sqrt{n}x,\sqrt{n}y )-\rho_\mu (t,x,y)\right|\leq  \textup{Err}(n;\mu,\epsilon,M,\tau).
\end{align*}

Furthermore,
\begin{align*}
\lim_{n\to\infty }\textup{Err}(n;\mu,\epsilon,M,\tau)=0.
\end{align*}

See the expression for $\rho_{\mu}$ in \eqref{Robin-kernel}.
\end{lem}

\begin{lem}\label{var_transition}
Fix $n\in\mathbb{N}$,\ $\tau\geq 1$, $\epsilon\in (0,1]$ and $K\geq 1$. Assume that $\Var [X_i] \leq Kn^{-\epsilon}  $ and that $\E[X] = \gamma = 1-  {\mu}/{\sqrt{n}}$. There exists $c(n;\tau,\epsilon,K)$ such that the following statement holds. For any integers $ m\in [1,n\tau]_\mathbb{Z} $ and $(x, y)\in\mathbb{N}^2_0 $, we have
\[
\Var [p_X(0,m,x,y)]\leq    c(n;\tau,\epsilon,K)    p_{\max\{1,\gamma^2\}}^2(0,m,x,y) .
\]
Furthermore,
$$\lim_{n\to\infty}c(n,\tau,\epsilon,K)=0.$$
\end{lem}

\section{Proof of the Main theorem}\label{proof-section}
For simplicity, we first treat the case $t=1$ and explain how to proceed with general $t>0$. In order to prove the convergence of the point-to-point partition functions $Z^{\omega,X}(n,[x\sqrt{n}]_n;\beta n^{-1/4})$ as in \eqref{p2p-part-function}, we begin with identifying $Z ^{\omega,X}(n,[x\sqrt{n}]_n;\beta n^{-1/4})$ with a U-statistics of $p_{X,k}(n,[x\sqrt{n}]_n;\bi,\mathbf{y} )$ as in \eqref{k-fold-prob} and then use the techniques of U-statistics. 

As $p_{X,k }(\floor*{nt},[x\sqrt{n}]_{\floor*{nt}};\bi,\mathbf{y})$ is only defined on lattice points $(\bi,\mathbf{y}) \in D_k^{\floor*{nt}}\times \N_0^k$, which verify the parity condition, we will interpolate the discrete transition kernel $p_{X,k }(\floor*{nt},[x\sqrt{n}]_{\floor*{nt}};\bi,\mathbf{y})$ to be a $L^2$ function on $[0,t]^k \times \R_{\geq 0} ^k$. Given $x \in \R_{\geq 0}$ and $i \in \N$, recall that $[x]_i$ defined in \eqref{def:parity} is the largest integer among the ones that are smaller than $x$ and are of the same parity as $i$. For a point $\bx \in \R_{\geq 0} ^k $ and $\bi \in D^{\floor*{nt}}_k$, define $[\bx] _{\bi} \in \mathbb{N}_{0} ^k$ by $([\bx ]_{\bi})_k = [x_k]_{i_k}$. 

Given $  (t,x) \in \mathbb{R}_{>0}\times\mathbb{R}_{\geq 0}$ and $(\bt,\bx)\in [0,t]^k  \times \mathbb{R}^k_{\geq 0}$, let $m=\floor*{nt}$, $y=[ \sqrt{n} x ]_m$, $\bi=\lfloor n \bt\rfloor$ and $\mathbf{y}=[\sqrt{n} \bx]_{\bi}$. Define the scaled extension $\nu_{X,k}^{n}$ as
\begin{equation}\label{def:k-fold}
\nu_{X,k}^{n}(t,x; \bt,\bx)\coloneqq  2^{-(k+1)}n^{(k+1)/2}p_{X,k }(m,y; \bi,\mathbf{y})\cdot{\mathbbm{1}}\{ \bi \in D_k^m\}.
\end{equation}

Note that now $\nu_{X,k}^{n} $ also takes care of the diffusive scaling for $p_{X,k} $. Under above definitions, $\nu_{X,k}^{n}$ is constant on the rectangles of $\mathcal{R}_k^n$. Note that for $\bi \in E_k^n, \mathbf{y}  \in \mathbb{N}_{0} ^k$ such that $\bi \leftrightarrow \mathbf{y}$,
\begin{align*}
\nu_{X,k}^{n}(t,x; n^{-1} \bi ,  n^{-1/2} \mathbf{y} )=2^{-(k+1)}n^{(k+1)/2} p_{X,k }(m,y; \bi,\mathbf{y} )\cdot \mathbbm{1}\{\bi \in D_k^m\}. 
\end{align*} 

Recall the definition of $\mathcal{S}_k^n$ as in \eqref{Def-U-stats} and note that $\nu_{X,k}^{n}$ is constant on the rectangles of $\mathcal{R}_k^n$ and zero elsewhere, we compute the U-statistics of $\nu_{X,k}^{n}(1,x;\cdot,\cdot)$ as follows,
\begin{align*}
\mathcal{S}_k^n(\nu_{X,k}^{n}(1,x;\cdot,\cdot) ;\omega)&= 2^{k/2 }  \sum_{\bi \in E_k^n} \sum_{\mathbf{y}  \in \mathbb{N}_{0}^k} \nu_{X,k}^{n} \left(1,x;n^{-1} \bi ,n^{-1/2}  \mathbf{y}  \right)\cdot \omega(\bi, \mathbf{y}) \cdot \mathbbm{1}\{\bi \leftrightarrow \mathbf{y}\} \\
 &= 2^{-k/2-1}n^{(k+1)/2} \sum_{\bi \in D_k^n} \sum_{\mathbf{y} \in \mathbb{N}_{ 0}^k} p_{X,k }(n,y;\bi,\mathbf{y})\cdot \omega(\bi, \mathbf{y}).
\end{align*}
Here $y=[x]_n$ and the parity condition is handled by the $p_{X,k}$ and summation is over $\bi \in D^n_k$. We could rewrite the modified point-to-point partition function as
\begin{align}\label{p2p_Ustat}
\mathfrak{Z}^{\omega,\gamma} (n,[x\sqrt{n}]_n;\beta n^{-1/4})= {2}n^{-1/2} \sum_{k=0}^n 2^{ k /2}\beta^k n^{-3k/4} \mathcal{S}_k^n( \nu_{X,k}^{n}(1,x;\cdot,\cdot);\omega).
\end{align}

The following two lemmas seek to bound $\nu_{X,k}^{n}$. Lemma~\ref{p_k^n_L_2} gives the $L^2$ bound and $L^2$ convergence of $\nu_{\gamma ,k}^{n}(t,x;\bt,\bx)$. 
\begin{lem}
Fix  $\mu\in\mathbb{R}$ and $\tau\geq 1$. There exist a constant $B_1(\mu,\tau)$ and a universal constant constant $C_1$ such that the following statement holds. For any $n\in\mathbb{N}$, define 
\begin{equation*}
\Theta_n(t_1,t_2,x_1,x_2)\coloneqq B_1\max\{t_2-t_1,2n^{-1}\}^{-1/2}e^{-(x_2-x_1)^2/(C_1\max\{t_2-t_1,2n^{-1}\})}. 
\end{equation*}
Let $\gamma=1-\mu/\sqrt{n}$. Then for all $n\geq 1$, $t\in (0,\tau]$, $ x\in\mathbb{R}$ and $(\bt,\bx)\in\Delta_k(t)\times\mathbb{R}^k$, we have
\begin{equation*}
\nu_{\gamma ,k}^{n}(t,x;\bt,\bx) \leq \mathbf{F}_k[\Theta_n](t,x;\bt,\bx).
\end{equation*}
Here the k-fold operator $\mathbf{F}_k$ is define in \eqref{def:Fk}.
\end{lem}

\begin{proof}
Let $m=\floor*{nt}$ $y= [\sqrt{n} x]_m$, $\bi=\floor* {n\bt}$ and $\mathbf{y}= [\sqrt{n} \bx]_{\bi}$. Without loss of generality we may assume $\bi\in D^m_k$ as otherwise $\nu_{\gamma ,k}^{n}(t,x;\bt,\bx)=0$. In particular, $m-i_k\geq 1$ and $i_j-i_{j-1}\geq 1$. By the definition of $\nu_{\gamma ,k}^{n}(t,x;\bt,\bx)$, it suffices to show that
\begin{align*}
2^{-1}n^{1/2}p_{\gamma}(i_{j-1},i_j,y_{j-1},y_j)\leq &\Theta_n(t_{j-1},t_j,x_{j-1},x_j),\\
2^{-1}n^{1/2}p_{\gamma}(i_{k},m,y_{k},y )\leq &\Theta_n(t_{k},t,x_{k},x).
\end{align*} 

We give the proof for the first inequality. The proof for the second is identical. From Lemma \ref{transition_p_unform_upperbd}, 
\begin{align*}
2^{-1}n^{1/2}p_{\gamma}(i_{j-1},i_j,y_{j-1},y_j)\leq Bn^{1/2}   (i_j-i_{j-1}) ^{-1/2}e^{-(y_j-y_{j-1})^2/[C(i_j-i_{j-1})]}.
\end{align*}

We assume first that $t_j-t_{j-1}\geq 2n^{-1}$. Then  
$$ (t_j-t_{j-1})/2\leq n^{-1}(i_j-i_{j-1})\leq 2(t_j-t_{j-1}).$$

Together with 
$$n(x_j-x_{j-1})^2\leq 2(y_{j}-y_{j-1})^2+4,$$
The assertion follows. The proof for $t_j-t_{j-1}< 2n^{-1}$ is similar by using 
$$n^{-1}\leq  n^{-1}(i_j-i_{j-1})\leq 3n^{-1}.$$

The proof is finished.
\end{proof}

\begin{lem}\label{p_k^n_L_2}
Fix $\mu\in\mathbb{R}$ and $\tau\geq 1$. There exists a constant $B_2(\mu,\tau)$ such the the following statement holds.  For all $n\in \mathbb{N}$, let $\gamma=1-\mu/\sqrt{n}$. For all $k\geq 1$, $t\in (0,\tau]$ and $x\in \R_{\geq 0}$, we have
\begin{align}
\label{p_k^n_L_2_part1} & \norm{ \nu_{\gamma ,k}^{n}(t,x;\cdot,\cdot) }_{L^2}^2 \leq t^{k/2-1} e^{-x^2/  [C_1 \max\{t,2k/n\}]  } B_2(\mu,\tau)^{k}/\Gamma((k+1)/2)  ,\\
\label{p_k^n_L_2_part2} &\lim_{n\rightarrow \infty}\norm{ \nu_{\gamma ,k}^{n}(t,x;\cdot,\cdot) - \rho_{\mu,k}(t,x;\cdot,\cdot)}_{L^2}=0.
\end{align}
\end{lem}
\begin{proof}
We start with \eqref{p_k^n_L_2_part1}. By a direct computation,
\begin{align*}
\mathbf{F}_k[\Theta_n](t,x;\bt,\bx)^2= B_1^{k+1}(\max\{t-t_{k},2n^{-1}\})^{-1/2}\prod_{j=1}^{k}(\max\{t_j-t_{j-1},2n^{-1}\})^{-1/2}\\
\times  \mathbf{F}_k[\Theta_n](1,\sqrt{2} x;\bt,\sqrt{2}\bx). 
\end{align*}
Through change of variables, for any $\bt\in \Delta_k(t)$,
\begin{align*}
\int_{\mathbb{R}^k}\mathbf{F}_k[\Theta_n](t,\sqrt{2} x;\bt,\sqrt{2}\bx)d\bx\leq  B^{k+1} t^{-1/2}e^{-x^2/[C_1 \max\{t,2k/n\}]}.
\end{align*}

For simplicity, we denote $\bar{t}=\max\{t,2k/n\}$. Thus
\begin{align*}
\sup_{n\in\mathbb{N}} \norm{\nu_{\gamma ,k}^{n}(t,x;\cdot,\cdot) }_{L^2}^2\leq  &B^{k}t^{-1/2}e^{-x^2/(C_1 \bar{t})}\int_{\Delta_k(t)}(\max\{t-t_{k},2n^{-1}\})^{-1/2}\prod_{j=1}^{k}(\max\{t_j-t_{j-1},2n^{-1}\})^{-1/2} d\bt\\
\leq &B^{k}t^{-1/2}e^{- x^2/(C_1\bar{t} )}\int_{\Delta_k(t)}( t-t_{k} )^{-1/2}\prod_{j=1}^{k}( t_j-t_{j-1}  )^{-1/2} d\bt\\
\leq & t^{k/2-1} e^{- x^2/(C_1\bar{t} )} B^{k}/\Gamma((k+1)/2).
\end{align*}
Here we have used (see \cite[Section 3.4]{AKQ})
\begin{align*}
\int_{\Delta_k}( t-t_{k} )^{-1/2}\prod_{j=1}^{k}( t_j-t_{j-1}  )^{-1/2} d\bt =\pi^{(k+1)/2}/\Gamma((k+1)/2).
\end{align*}

Next, we turn to showing \eqref{p_k^n_L_2_part2}. By the local limit theorem Lemma \ref{pointwise_convergence_tran_barrier}, $\nu_{\gamma ,k}^{n}(t,x;\cdot,\cdot) $ converges to $\rho_{\mu,k}(t,x;\cdot,\cdot)$ pointwisely. By the argument above we see that $\mathbf{F}_k[\Theta_n](t,x;\cdot,\cdot)$ converges to $\mathbf{F}_k[\Theta_\infty](t,x;\cdot,\cdot)$ in $L^2$. Here
\begin{align*}
\Theta_\infty(t_1,t_2,x_1,x_2)\coloneqq B_1(t_2-t_1)^{-1/2}e^{-(x_2-x_1)^2/(C_1(t_2-t_1) ) }. 
\end{align*}

Thus \eqref{p_k^n_L_2_part2} follows by the dominated convergence theorem.
\end{proof}

By identifying $\mathfrak{Z} ^{\omega ,\gamma}(n,[x\sqrt{n}]_n;\beta n^{-1/4})$ with the U-statistics as in \eqref{p2p_Ustat}, we are ready to prove the main Theorem~\ref{main_thm} in a few steps as follows. 

\begin{proof}[Proof of Theorem~\ref{main_thm}]
Define the environment field $\omega_n$ by
\begin{align}\label{tilde_omega}
e^{\beta n^{-1/4} \omega(i,x) - \lambda(\beta n^{-1/4})} = 1 + \beta
n^{-{1}/4}  {\omega}_n(i,x).
\end{align}
Note that as $\E[e^{\beta_0\omega}]<\infty$, $\lambda(\beta n^{-1/2})$ is well-defined as $\beta n^{-1/4}\leq \beta_0$. From the definition of $\lambda(\beta n^{-1/4})$, we have $\E [{\omega}_n] = 0$. It is straightforward to check that $\E[ {\omega}_n^2]=1+O(n^{-1/4})$. Hence $\omega_n$ satisfies \eqref{omega}. Moreover we have
\begin{equation*}
\begin{split}
&2^{-1}n^{1/2} e^{ - n\lambda(\beta n^{-1/4})} Z ^{\omega,X} \left(n,[x\sqrt{n}]_n;\beta n^{-{1}/{4}}  \right)\\
 = &2^{-1}n^{1/2} \E_R \left[
\prod_{i=0}^n \left(1 + \beta n^{-{1}/{4}} {\omega}_n(i,S_i) \right)\mathbbm{1}\{S_n = [x\sqrt{n}]_n\} \right]\\
 =&2^{-1}n^{1/2} \mathfrak{Z} ^{ {\omega}_n,X} \left(n,[x\sqrt{n}]_n ;\beta n^{-{1}/{4}} \right).
 \end{split}
\end{equation*}

Step 1: Fix $x\in\mathbb{R}_{\geq 0}$. We first prove the convergence of $2^{-1}n^{1/2} \mathfrak{Z} ^{\omega_n,\gamma}(n,[x\sqrt{n}]_n;\beta n^{-1/4})$. By \eqref{p_k^n_L_2_part1} and \eqref{p_k^n_L_2_part2},
\begin{align*}
\|\rho_{\mu,k}(1,x;\cdot,\cdot) \|^2_{L^2}\leq e^{-2x^2/C_1}B_2(\mu,1)^2/\Gamma((k+1)/2).
\end{align*}

It is easy to see that \eqref{GG} holds. Hence by Lemma \ref{lemma_U-stat} it follows that for all $\beta >0$, as $n\rightarrow \infty$,
\begin{align}\label{rho_series}
\sum_{k=0}^{\infty} 2^{k/2}\beta^k
n^{-3k/4} \mathcal{S}_k^n(\rho_{\mu,k}(1,x;\cdot,\cdot);\omega_n)
\mathop{\longrightarrow}^{(d)} z_{\sqrt{2}\beta}(1,x).
\end{align}
See the chaos expansion of $z_{\sqrt{2}\beta}(1,x)$ in \eqref{mild_solution} . Now it suffices to show that the difference
\begin{align*}
J&\coloneqq\sum_{k=0}^{\infty} 2^{k/2}\beta^k
n^{-3k/4} \mathcal{S}_k^n(\rho_{\mu,k}(1,x;\cdot,\cdot) ;\omega_n)-2^{-1}n^{1/2}\mathfrak{Z} ^{\omega_n,\gamma}(n,[x\sqrt{n}]_n;\beta n^{-{1}/{4}}).
\end{align*}
converges to $0$ in $L^2$. By splitting the above series and applying linearity of $\mathcal{S}_k^n$, we have
\begin{align*}
J &=\sum_{k=0}^{\infty} 2^{k/2}\beta^k
n^{-3k/4} \mathcal{S}_k^n(\rho_{\mu,k}(1,x;\cdot,\cdot) ;\omega_n)- \sum_{k=0}^n 2^{k/2}\beta^k
n^{-3k/4} \mathcal{S}_k^n(  \nu_{\gamma ,k}^{n}(1,x;\cdot,\cdot) ;\omega_n) \\
&=\sum_{k=0}^n 2^{k/2}
\beta^k n^{-3k/4} \mathcal{S}_k^n \bigl(
\rho_{\mu,k}(1,x;\cdot,\cdot)  -   \nu_{\gamma ,k}^{n}(1,x;\cdot,\cdot) ;\omega_n \bigr) +
\sum_{k=n+1}^{\infty} 2^{k/2} \beta^k
n^{-3k/4} \mathcal{S} _k^n(\rho_{\mu,k}(1,x;\cdot,\cdot) ;\omega_n).
\end{align*}
Because $\mathcal{S} _k^n(\rho_{\mu,k}(1,x;\cdot,\cdot) ;\omega_n)$ are independent for different $k$, by Lemma~\ref{lemma-U-bounds} the second moment of the second term is bounded from above by
\[
\sum_{k = n+1}^{\infty} \E[\omega^2_n]^k2^k  \beta^{2k} \Vert \rho_{\mu,k}(1,x;\cdot,\cdot) \Vert ^2_{L^{2} }\leq \sum_{k = n+1}^{\infty} \E[\omega^2_n]^k2^k \beta^{2k}  e^{-2x^2/C_1}B_2(\mu)^k/\Gamma((k+1)/2).
\]
Thus the second term converges to zero as $n$ goes to infinity. We now turn to the first term. By Lemma~\ref{lemma-U-bounds} we have
\[ 
\mathbb{E} \left[\left(\sum_{k=0}^n 2^{k/2}
\beta^k n^{-3k/4} \mathcal{S}_k^n \bigl(
\rho_{\mu,k}(1,x;\cdot,\cdot) - \nu_{\gamma ,k}^{n}(1,x;\cdot,\cdot) \bigr)\right)^2\right]
\leq
\sum_{k=0}^n \E[\omega^2_n]^k 2^{k}
\beta^{2k} \bigl\Vert \rho_{\mu,k}(1,x;\cdot,\cdot) -  
\nu_{\gamma ,k}^{n}(1,x;\cdot,\cdot) \bigr\Vert _{L^2}^2.
\]
Lemma~\ref{p_k^n_L_2} shows that for any $k$, as $n \rightarrow \infty$, $\bigl\Vert \rho_{\mu,k}(1,x;\cdot,\cdot) -
  \nu_{\gamma ,k}^{n}(1,x;\cdot,\cdot) \bigr\Vert _{L^2}^2 \rightarrow 0$. Together with $\bigl\Vert  \rho_{\mu,k}(1,x;\cdot,\cdot) -
  \nu_{\gamma ,k}^{n}(1,x;\cdot,\cdot) \bigr\Vert _{L^2}^2\leq 4e^{- x^2/C_1\max\{1,2k/n \}}B_2(\mu,1)^k/\Gamma((k+1)/2)$ from Lemma \ref{p_k^n_L_2_part1}, it follows by dominated convergence theorem that
\[
\lim_{n \to\infty} \sum_{k=0}^n
\E[\omega^2_n]^k 2^k \beta^{2k} \bigl\Vert \rho_{\mu,k}(1,x;\cdot,\cdot) -  \nu_{\gamma ,k}^{n}(1,x;\cdot,\cdot) \bigr\Vert _{L^2}^2
= 0.
\]
We then conclude that
\begin{align*}
2^{-1}n^{1/2} \mathfrak{Z}^{\omega_n,\gamma}(n,[x\sqrt{n}]_n;\beta n^{-1/4})\mathop{\longrightarrow}^{(d)} z_{\sqrt{2}\beta}(1,x).
\end{align*}

Step 1 is finished.\\

Step 2: We now turn to demonstrating convergence of $\mathfrak{Z}^{\omega_n,X}(n,[x\sqrt{n}]_n;\beta n^{-1/4})$ where randomness is also present at the boundary random environment. It suffices to show
\[
2^{-1}n^{1/2}\left( \mathfrak{Z} ^{\omega_n,X}(n,[x\sqrt{n}]_n;\beta n^{-1/4}) -   \mathfrak{Z} ^{\omega_n,\gamma}(n,[x\sqrt{n}]_n;\beta n^{-1/4}) \right) \xrightarrow{(d)} 0.
\]
We have 
\begin{align*}
&2^{-1}n^{1/2}\left( \mathfrak{Z}^{\omega_n,X}(n,[x\sqrt{n}]_n;\beta n^{-1/4}) -   \mathfrak{Z} ^{\omega_n,\gamma}(n,[x\sqrt{n}]_n;\beta n^{-1/4}) \right)\\
&=\sum_{k=0}^n 2^{k/2}
\beta^k n^{-3k/4} \mathcal{S}_k^n \bigl(  \nu_{X,k}^{n}(1,x;\cdot,\cdot)  - \nu_{\gamma ,k}^{n}(1,x;\cdot,\cdot)   ;\omega_n\bigr).
\end{align*}

By Lemma \ref{lemma-U-bounds},  
\begin{equation*}
\begin{split}
\Var\left[2^{-1}n^{1/2}\left( \mathfrak{Z}_n^{\omega,X}([x\sqrt{n}]_n;\beta n^{-1/4}) -  \mathfrak{Z}_n^{\omega,\gamma}([x\sqrt{n}]_n;\beta n^{-1/4}) \right) \right] \\
\leq   \sum_{k=0}^{n} \E[\omega_n^2]^k2^{k}\beta^{2k}\int_{\Delta_k \times\R^k_{\geq 0}}\E \left[   \nu_{X ,k}^{n}(1,x;\vt,\vx)  - \nu_{\gamma ,k}^{n}(1,x;\vt,\vx) \right]^2 d\vt d\vx.  
\end{split}
\end{equation*}

Recall the definition for $\nu^{n}_{X,k }(1,x;\bt,\bx)$ as in \eqref{def:k-fold} and the definition for the k-fold transition kernel $p_{X,k}$ as in \eqref{k-fold-prob}. Fix $n\in\mathbb{N}$. Let $y=[ \sqrt{n} x ]_n$, $\bi=\lfloor n \bt\rfloor$ and $\mathbf{y}=[\bx]_{\bi}$. Without loss of generality we may assume $\bi\in D^n_k$. As $\E \left[ \nu^{n}_{X,k }  \right]=    \nu^{n}_{\gamma,k }  $, it follows that
\begin{align*}
& \quad\;2^{ 2(k+1)}n^{k+1} \E \left[  (\nu^{n}_{X,k } - \nu^{n}_{\gamma,k })(1,x;\bt,\bx) \right]^2\\
&= \E[ p^2_{X}(i_k, n, y_{k} ,y ) ]\prod _{j=1}^k \E [ p^2_{X}(i_{j-1},i_j ,y_{j-1} , y_{j})]-p^2_{\gamma}(i_k, n, y_{k} ,y )\prod _{j=1}^k  p^2_{\gamma}(i_{j-1},i_j, y_{j-1} , y_{j} )
\end{align*}
By Lemma~\ref{var_transition}, under the assumption $\Var(X_i)\leq Kn^{-\epsilon}$, 
$ \Var[p_X(m,n,x,y)] = c(n ;\epsilon,K) p^2_{\max\{1,\gamma^2\}}(m,n,x,y)$ with $\lim_{n\to\infty} c(n; \epsilon,K)=0$. Hence
$$ 0\leq \E[  p^2_X(m,n,x,y)] - p^2_{\gamma}(m,n,x,y) \leq    c(n ;\epsilon,K) p^2_{\max\{1,\gamma^2\}}(m,n,x,y).$$ 

By taking $n$ large enough, we may assume $c(n;\epsilon,K)\leq 1$. Then

\begin{align*}
0\leq &\E[ p^2_{X}(i_k, n, y_{k} ,y ) ]\prod _{j=1}^k \E [ p^2_{\gamma}(i_{j-1},i_j, y_{j-1} , y_{j}]-p^2_{\gamma}(i_k, n, y_{k} ,y )\prod _{j=1}^k  p^2_{\gamma}(i_{j-1},i_j, y_{j-1} , y_{j} )\\
\leq &2^kc(n;\epsilon,K)  p^2_{\max\{1,\gamma^2\} }(i_k, n, y_{k} ,y )\prod _{j=1}^k  p^2_{\max\{1,\gamma^2\}}(i_{j-1},i_j, y_{j-1} , y_{j} )\\
=&2^kc(n;\epsilon,K)\nu_{\max\{1,\gamma^2\},k}^n(1,x,\bt,\bx).
\end{align*}

By \eqref{p_k^n_L_2_part1}, we deduce
\begin{align*}
\int_{[0,1]^k\times\R^k_{\geq 0}}\E \left[   (\nu^{n}_{X,k } - \nu^{n}_{\gamma,k })(1,x,\bt,\bx) \right]^2 d\vt d\vx\leq  c(n;\epsilon,K)  B^{k}/\Gamma((k+1)/2).
\end{align*}

Hence
\begin{align*}
&\Var\left[2^{-1}n^{1/2}\left( \mathfrak{Z} ^{\omega,X}(n,[x\sqrt{n}]_n;\beta n^{-1/4}) -   \mathfrak{Z} ^{\omega,\gamma}(n,[x\sqrt{n}]_n;\beta n^{-1/4}) \right) \right]\\
\leq& c(n,\epsilon,K)    \sum_{k=0}^\infty \E[\omega^2_n]^k B^{k}\beta^{2k}    /\Gamma((k+1)/2)\to 0.
\end{align*}

As a result, $$ 2^{-1}n^{1/2} e^{ -n \lambda(\beta n^{-1/4})} Z ^{\omega,X} \left(n,[x\sqrt{n}]_n;\beta n^{-{1}/{4}}  \right) = 2^{-1}n^{1/2} \mathfrak{Z} ^{\omega_n,X}(n,[x\sqrt{n}]_n;\beta n^{-1/4}) \xrightarrow{{\hspace*{0.2cm}(d)}\hspace*{0.2cm}}z_{\sqrt{2}\beta}(1,x).$$ 

This proves the one point convergence of $2^{-1}n^{1/2} e^{ -n \lambda(\beta n^{-1/4})} Z ^{\omega,X} \left(n,[x\sqrt{n}]_n;\beta n^{-{1}/{4}}  \right)$.\\

Note that for all $t>0$, Lemma \ref{var_transition} and Lemma \ref{p_k^n_L_2} hold. Hence the argument above actually yields the convergence of $2^{-1}n^{1/2} e^{ -nt \lambda(\beta n^{-1/4})} Z ^{\omega,X} \left(\floor*{nt} ,[x\sqrt{n}]_{\floor*{nt}};\beta n^{-{1}/{4}}  \right)$ to $z_{\sqrt{2}\beta}(t,x)$ for arbitrary $t>0$ and $x\in\mathbb{R}$. Furthermore, by the joint convergence in Lemma \ref{lemma_U-stat}, the finite dimensional convergence also follows.\\

Step 3: Now in order to show the weak convergence as a process, it suffices to show the tightness of the above process, which could be done by a similar argument as in \cite[Appendix B]{AKQ}. They first deduced an integral form in terms of the random walk transition kernel for the modified point-to-point partition function $\mathfrak{Z}(x,k)$ from the discrete stochastic heat equation that $\mathfrak{Z}(x,k)$ satisfies and then developed the modulus of continuity for the partition function with estimates for heat kernel. In our case, for deterministic $X_i\equiv \gamma$, we could derive a similar integral form for the point-to-point partition function but in terms of transition kernel for half-line random walk with a barrier at origin and then the similar estimates follow given that Robin heat kernel has similar decay behavior as standard heat kernel as in Lemma~\ref{transition_p_unform_upperbd}. For $X_i$ under the assumption of Theorem \ref{main_thm}, from Remark \ref{third_moment}, we have that $\mathbb{E}[|\nu^1_{X,1} -\nu^1_{\gamma,1}|^\alpha]$ converges to zero in $L^1([0,1]\times \mathbb{R}_{\geq 0})$ for any $1\leq \alpha<3$. Here $\nu^1_{X,1}  $ and $\nu^1_{\gamma,1} $ are interpolated (random) transition kernel as in \eqref{def:k-fold}. This allows us to adapt the proof in \cite[Appendix B]{AKQ} to the current setting.

\end{proof}

\section{Application to log-gamma polymer models}\label{application-section}
In this section we consider the half-space log-gamma polymer model, as introduced in \cite{Sep}. We apply the main Theorem~\ref{main_thm} to the log-gamma polymer point-to-point partition function. The log-gamma polymer models in dimension $1+1$ are of significant importance among polymer models in the sense that integral formulas are discovered and steepest descent analysis is allowed, see \cite{BBC}. 

We start with defining the log-gamma polymer model. We first follow notations used in the literature and then translate it to fit our setting for the general polymers. 
Consider a half-quadrant $V\coloneqq\{(i,j) |  i\geq j, i,j \in \mathbb{N}_0\}$. Assign a log-gamma random environment $Y\coloneqq \{Y_{i,j}, i\geq j\}$ on $V$ as follows.
\begin{align}\label{weights-log-gamma}
Y_{i,i} \sim \textup{Inv-Gamma} (\sqrt{n}+\mu+1/2), \quad Y_{i,j} \sim \textup{Inv-Gamma} (2\sqrt{n}),\  \textup{for}\ i>j.
\end{align}
Here $\textup{Inv-Gamma}(\alpha)$ is the inverse gamma distribution with shape parameter $\alpha$ and scale parameter 1, and with density
\[
\frac{1}{\Gamma(\alpha)}x^{-\alpha-1}e^{-1/x}.
\]
These choices of parameters correspond to the diffusive scaling and critical scaling at the origin of the general polymers.

For an endpoint $(m,n) \in V$, define the point-to-point partition function by
\[
Z^Y_{m,n} \coloneqq \sum_{\substack{S:(0,0)\rightarrow (m,n)}}\prod_{(i,j)\in S} Y_{i,j},
\]
where we sum over the up-right paths $S$ from $(0,0)$ to $(m,n)$ which always stay in the half-quadrant $V$. Note that the probabilities of these paths do not sum to one since those paths having crossed boundary $x=y$ are not counted. 

To match with the general environment setting in the half-space regime with a barrier at origin, we need to rewrite the partition function $Z^Y$ in the same form as \eqref{modified-p2p}, i.e. expectation with respect to a reflected random walk measure. By taking
\[
\tilde{Y}_{i,i} = \frac{1}{2} Y_{i,i}, \quad\tilde{Y}_{i,j} =  Y_{i,j}, i>j,
\] 
we have
\begin{align}\label{p2p_log_g}
Z^Y_{m,n} &= 2^{m+n}\sum_{S:(0,0)\rightarrow (m,n)} 2^{-(m+n)}\cdot 2^{\#_S}\cdot2^{-\#_S}\cdot\prod_{(i,j)\in S} Y_{i,j} \\\nonumber
&= 2^{m+n}\E_R\left[\prod_{(i,j)\in S} \tilde{Y}_{i,j}\cdot \mathbbm{1}\{S(m+n)=(m,n)\} \right],\nonumber
\end{align}
where $\#_S$ is the number of times that path $S $ visits the boundary and $\E_R$ is the expectation with respected to the reflected random walk measure. 

Once again we omit the floor function when it does not cause ambiguity, e.g. $\floor*{nt},  [x\sqrt{n}]_{\floor*{nt}}$. The following convergence result holds for log-gamma polymers.
\begin{thm}\label{log_gamma_part_thm}
Let $Y_{2,1}, $ be a random variable distributed as in ~\eqref{weights-log-gamma}. The following convergence results hold for the half-space log-gamma polymer model as $n \rightarrow \infty$,
\begin{equation*}
(2^{-1}n^{1/2}) 2^{-\floor{nt}}\E[Y_{2,1}]^{-\floor{nt}}  \cdot Z_{\floor{ ( nt+x\sqrt{n})/2}, \floor{ (nt-x\sqrt{n})/2}} \xrightarrow[]{(d)}\;  z_{1}(t,x).
\end{equation*}
\end{thm}

\begin{proof}
From \eqref{p2p_log_g}, we have
\begin{align*}
&\quad 2^{-nt}\E[Y_{2,1}]^{-nt} \cdot Z_{\floor{ ( nt+x\sqrt{n})/2}, \floor{ (nt-x\sqrt{n})/2}} \\
&= \E_R\left[\left(\prod_{(i,j)\in S} {\E[Y_{2,1}]}^{-1}{\tilde{Y}_{i,j}}\right) \mathbbm{1}\left\{S(nt)=\left(\floor{ (nt+x\sqrt{n})/2}, \floor{ (nt-x\sqrt{n})/2}\right)\right\}\right].
\end{align*}
Define $\omega_n(i,j)$ for $i\geq j$ via
\begin{align*}
{\E[Y_{2,1}]}^{-1} {\tilde{Y}_{i,j}}&=:1+2^{-1/2}  {n^{-1/4}} \omega_n(i,j), i >j;\\
{\E[Y_{2,1}]}^{-1} {\tilde{Y}_{i,i}}&=:\gamma_n\left(1+2^{-1/2} n^{-1/4}  \omega_n(i,i)\right),
\end{align*}
where $\gamma_n \coloneqq 2^{-1}\E[Y_{i,i}]/\E[Y_{2,1}].$

In these notations, it follows that
\[
(2^{-1}n^{1/2}) 2^{-nt}\E[Y_{2,1}]^{-nt}  \cdot Z_{\floor{\frac{1}{2}( nt+x\sqrt{n})}, \floor{\frac{1}{2}(nt-x\sqrt{n})}} = 2^{-1}n^{1/2} \mathfrak{Z} ^{\omega_n,\gamma_n}(\floor*{nt} ,\floor*{\sqrt{n}x }; 2^{-1/2}  n^{-{1}/{4}} ).
\] 
The shows that the log-gamma partition function is equivalent to the scaled modified point-to-point partition function as in \eqref{modified-p2p} with $\beta = \frac{1}{\sqrt{2}}$.

Furthermore, it's clear that for $i\geq j$, $\E[\omega_n(i,j)] = 0$. And since 
$$\E[\textup{Inv-Gamma}(\alpha)] =  (\alpha -1)^{-1}, \quad \Var[\textup{Inv-Gamma}(\alpha)] =   (\alpha-1)^{-2}(\alpha-2)^{-1}, $$
we deduce,
\begin{align*}
\Var[\omega_n(i,j) ] &= 2n^{1/2} \Var[Y_{i,j}]/\E[Y_{i,j}]^2  = 1+O\left(n^{-1/2} \right), i>j;\\
\Var[\omega_n(i,i)] &= 2n^{1/2} {\Var[Y_{i,i}]}/{\E[Y_{i,i}]^2}= 2+O\left(n^{-1/2}  \right).\\
\gamma_n &  = 1- {\mu}/{\sqrt{n}}+O\left(n^{-1}\right).
\end{align*}

Note that now the weights $\omega_{i,j}$ on the off-diagonals are i.i.d. with mean zero and variance asymptotically one, the weights $\omega_{i,i}$ on the diagonal are also i.i.d. with mean zero but with variance asymptotically two. Also for $\gamma_n =1- {\mu}/{\sqrt{n}} +O\left(n^{-1} \right)$, we have the same local limit theorem as in Theorem~\ref{pointwise_convergence_tran_barrier}.

The rest of this proof follows exactly the same argument as in Theorem~\ref{main_thm}, with the role of U-statistics Lemma~\ref{lemma_U-stat} being replaced by Lemma~\ref{sum_sk_tilde_omega}. Hence the desired convergence for log-gamma polymer model holds. 
\end{proof}

\begin{appendix}
\section{Proofs for Section \ref{heat_kernel_estimates}}
In this section we prove the three lemmas in Section \ref{heat_kernel_estimates}, i.e. Lemmas~\ref{transition_p_unform_upperbd}, \ref{pointwise_convergence_tran_barrier} and \ref{var_transition}. The proofs rely on a few lemmas on estimates for random walks. The reader may skip these lemmas first and proceed directly to the proofs of Lemmas~\ref{transition_p_unform_upperbd}, \ref{pointwise_convergence_tran_barrier} and \ref{var_transition}. It will be further explained in the proofs that which lemmas will be applied.

Recall that $\gamma$ is the reflection rate, when $\gamma \leq 1$, $p_{\gamma}(m,m+n,x,y)\leq p(m,m+n,x,y)$, i.e the totally reflecting case, but when $\gamma >1$ the system will have mass coming in. Therefore we need to estimate how frequently the walker goes to the barrier in order to estimate the discrete transition kernel. 

Recall that the transition kernel $p_\gamma$ is defined as

\begin{align}\label{eq_p} 
p_{\gamma}(m,n,x,y)=\sum_{j=0}^{n-m}\gamma^j\mathbb{P}^{m,x}_R(N_{m,n}=j,S_n=y).
\end{align}
Here $N_{m,n}$ is the total visits to the origin as
\begin{align*} 
N_{m,n}(S)= \#\{i\in [m,n-1]_\mathbb{Z}\ |\ S_i=0\}.
\end{align*}
For the case $m=0$, we denote $N_{0,n}$ as $N_n$ to simplify the notation. 
The explicit form of $\mathbb{P}^{ m,x }_R(N_{m,n}=j,S_{n}=y)$, (see Lemma~\ref{tran_prob_barrier}), can be found in \cite[(27)]{Goo}.  We give a proof in the appendix for the reader's convenience. For $(n,z)\in\mathbb{N}_0\times \mathbb{Z}$, let $T(n,z)$ be the probability that a simple random walk on $\Z$ arrives at $x = z$ after $n$ jumps starting at origin. In other words,
\[
 T(n,z)\coloneqq \mathbb{P}^{0,0}(S_n=z). 
\]
\begin{lem}\label{tran_prob_barrier}
For any $(m,n,x,y)\in \mathbb{N}_0\times \mathbb{N}\times \mathbb{N}_0\times \mathbb{N}$,
\begin{align*}
\mathbb{P}^{ m,x }_R(N_{m,m+n}=j,S_{m+n}=y) =\left\{ \begin{array}{cc}
 T(n,y-x)-T(n,y+x) & j=0,\\
 \frac{2(y+x+j-1)}{n-j+1}T(n-j+1,y+x+j-1) & j\geq 1.
\end{array} \right.
\end{align*}
And
\begin{align*}
\mathbb{P}^{ m,x }_R(N_{m,m+n}=j,S_{m+n}=0) =\left\{ \begin{array}{cc}
 ( T(n-1,1-x)-T(n-1,1+x))/2 & j=0,\\
 \frac{  x+j }{n-j }T(n-j , x+j ) & j\geq 1.
\end{array} \right.
\end{align*}
Note that the expression takes different form for $y=0$ and $y\neq 0$.
\end{lem}
The following two lemmas provide bounds on $T(n,z)$ and follow from computations through Stirling formula. The author did not find a reference for such results so proofs are provided in the next section.
\begin{lem}\label{CLT}
There exists a universal constant $C_2>0$ such that the following statement holds. For any $n\in\mathbb{N}$, $z\in\mathbb{Z}$, $z\leftrightarrow n$ and $|z|\leq n$, let $E(n,z):= {|z|^3}/{n^2}+ {1}/{n}$. Then
\begin{equation}\label{CLT_2}
e^{-C_2E(n,z)} \leq  2^{-1}(2\pi n)^{1/2}e^{z^2/(2n)}  T(n,z)\leq e^{ C_2E(n,z)}.
\end{equation}
\end{lem}

\begin{lem}\label{binominal_upper_bound}
There exists a universal constant $C_3>0$ such that the following statement holds. For any $n\in\mathbb{N}$, $z\in\mathbb{Z}$ and $z\leftrightarrow n$, we have 
\begin{equation}\label{upper_bound}
T(n,z)\leq C_3n^{-1/2} e^{ -z^2/( C_3n) }.
\end{equation}
\end{lem}

The following Lemma~\ref{1D_local_time_conditioned} and Lemma~\ref{pt_to_pt} seek bound for the expression in \eqref{eq_p}.
\begin{lem}
There exists a universal constant $C_4>0$ such that the following statement holds. For any $n\geq 1$, $x,y\in \mathbb{N}_0$ with $x-y\leftrightarrow n$ and $k\geq 0$, we have
\begin{align}\label{1D_local_time_conditioned}
\mathbb{P}_R^{0,x}(N_n\geq k|S_n=y)\leq C_4  e^{-k^2/(C_4n)}.
\end{align}
\end{lem}

\begin{proof} We first consider the case that $n$ is even and $x=y=0$. From Lemma \ref{tran_prob_barrier}, for any $k\geq 1$,
\begin{align*}
\mathbb{P}_R^{0,0}(N_n\geq k,S_n=0)=&\sum_{j\geq k}^{n/2}\frac{j}{n-j}T(n-j,j) \leq  \sum_{j\geq k}^{n/2}2C_3 (j/n)n^{-1/2}e^{-j^2/(2C_3n)}  \\
                                   = & 2C_3n^{-1/2} \sum_{j\geq k}^{n/2} (j/\sqrt{n}) e^{-(j/\sqrt{n})^2/(2C_3)} \cdot n^{-1/2},
\end{align*}
where the inequality follows from Lemma \ref{binominal_upper_bound}. 

Let $M_0>0$ be the number such that the function $se^{-s^2/(2C_3)}$ is decreasing for $s\geq M_0$. If $k< M_0\sqrt{n}$,  \eqref{1D_local_time_conditioned} holds easily as the right hand side can be made larger than $1$ with suitable $C_4$. 

Now we may assume $k\geq M_0\sqrt{n} $. By the integral test,
\begin{align*}
\mathbb{P}_R^{0,0}(N_n\geq k,S_n=0) \leq & C_3n^{-1/2} \int^\infty_{k/\sqrt{n}}se^{-s^2/(2C_3)} ds                                  \leq    Cn^{-1/2} e^{-k^2/(Cn)}.
\end{align*}
From \eqref{CLT_2}, $\mathbb{P}_R^{0,0}(S_n=0)=T(n,0)\geq 2(2\pi n)^{-1/2}e^{-C_2/n}$. Hence
\begin{align*}
\mathbb{P}_R^{0,0}(N_n\geq k | S_n=0)= & {\mathbb{P}_R^{0,0}(N_n\geq k,S_n=0)}\big/{\mathbb{P}_R^{0,0}(S_n=0)} \leq C e^{-k^2/(Cn)} .
\end{align*}
Thus \eqref{1D_local_time_conditioned} follows. 

Next, we consider general $x,y$ and $n$. Conditioning on the first and the last time the random walk bridge returns to the origin, we have for any $k\geq 2$,
\begin{align*}
\mathbb{P}_R^{0,x}(N_n\geq k | S_n=y) \leq &\max_{1\leq j\leq n } \mathbb{P}_R^{0,0}(N_{ j} \geq k-1 | S_{ j}=0).
\end{align*}
The change from $k$ to $k-1$ is necessary as $N_j$ ignores the zero at the end. Then \eqref{1D_local_time_conditioned} follows by the previous special case $x=y=0$.
\end{proof}

\begin{lem}\label{pt_to_pt}
For any $\mu\in\mathbb{R}$ and $\tau>0$, there exist  a constant  $ B_3(\mu,\tau)$ and a universal constant $C_5$ such that the following statement holds. For any $M\geq 0$, $n\in\mathbb{N}$, $m\in [1,\tau n]_\mathbb{Z}$ and $(x, y)\in \mathbb{N}^2_0$, let $\gamma=1-\mu/\sqrt{n}$. Then 
\begin{align*}
\sum_{k\geq M\sqrt{n}}\gamma^k \mathbb{P}^{0,x}_R (N_m(S)=k|S_m=y)\leq B_3(\mu,\tau)e^{- nM^2/(C_5 m )}.
\end{align*}
\end{lem}

\begin{proof} 
As $\gamma$ is decreasing in $\mu$, we can without loss of generality assume that $\mu\leq 0$. 
By \eqref{1D_local_time_conditioned} and  $\gamma \leq e^{|\mu|/\sqrt{n}}$, we obtain
\begin{align*}
&\sum_{k\geq M\sqrt{n}} \gamma^{k} \mathbb{P}^{0,x}_R (N_m=k|S_m=y)\\
=& (1-\gamma^{-1})\sum_{k\geq M\sqrt{n}+1} \gamma^{k} \mathbb{P}^{0,x}_R (N_m\geq k|S_m=y)+ \gamma^{M\sqrt{n}} \mathbb{P}^{0,x}_R (N_m\geq M\sqrt{n} |S_m=y)\\
\leq & C_4|\mu|n^{-1/2} \sum_{k\geq M\sqrt{n}} e^{-k^2/(C_4m)+k|\mu|/\sqrt{n}}   + C_4e^{-nM^2/(C_4m)+M|\mu|}.
\end{align*}

Here we have used summation by parts. As $m\leq \tau n$, $ k|\mu|/\sqrt{n}  \leq  k^2/(2C_4m)+    \tau C_4|\mu|^2/2$ and $ M|\mu|\leq  nM^2/(2C_4m) +  \tau C_4|\mu|^2/2$. Hence the above is bounded by
\begin{align*}
C_4 e^{\tau C_4|\mu|^2/2}\left(e^{-nM^2/(2C_4m)}+ |\mu|n^{-1/2}\sum_{k\geq M\sqrt{n}} e^{-k^2/(2C_4m)} \right).
\end{align*}

By the integral test,
\begin{align*}
n^{-1/2}\sum_{k\geq M\sqrt{n}} e^{-k^2/(2C_4m)}\leq (m/n)^{1/2}\int_{(n/m)^{1/2}M}^\infty e^{-s^2/(2C_4)}ds\leq C\tau^{1/2} e^{-nM^2/(2C_4m)}.
\end{align*}

Thus the assertion follows by putting the above together.
\end{proof}
\begin{proof}[proof of Lemma \ref{transition_p_unform_upperbd}]
By taking $M=0$ in Lemma \ref{pt_to_pt}, 
\begin{equation*}
 p_\gamma(0,m,x,y)=\sum_{k\geq 0}\gamma^k \mathbb{P}^{0,x}_R (N_m(S)=k|S_m=y)\mathbb{P}^{0,x}_R ( S_m=y)  \leq   B_3(\mu,\tau)\mathbb{P}^{0,x}_R ( S_m=y).
\end{equation*}
Together with Lemma \ref{binominal_upper_bound} and
\begin{align*}
\mathbb{P}^{0,x}_R ( S_m=y)=\left\{ \begin{array}{cc}
T(m,y-x)+T(m,y+x) & y\neq 0,\\
T(m,x)            & y=0.
\end{array} \right.
\end{align*}
The upper bound for $p_\gamma(0,m,x,y)$ follows.
\end{proof}

We are ready to prove the local limit theorem for $p_{\gamma}(m,m+n,x,y)$. Note that $p_{\gamma}(m,m+n,x,y)$ is indeed time-homogeneous and we may without loss of generality assume $m=0$.

\begin{proof}[proof of Lemma \ref{pointwise_convergence_tran_barrier}]
To simplify the notation, we adapt the convention that $C$ represents universal constants and $B$ represents constants that depend on $\mu,\epsilon,\tau$ and $M$. We adapt the notation that $A_1=A_2e^{\pm A_3}$ stands for $A_2 e^{-A_3}\leq  A_1\leq  A_2e^{ A_3}$. In particular, we can rewrite \eqref{CLT_2} as 
\begin{equation}\label{LLT}
T(n,z)=2(2\pi n)^{-1/2}e^{-z^2/(2n)\pm C_2E(m,z)}.
\end{equation}

We focus on the case that $y\neq 0$. The proof for $y=0$ is similar. Furthermore, we assume $n\geq n_0$ with $n_0$ large enough depending on $\mu, \epsilon,\tau$ and $M$. The exact value of $n_0$ may increase from line to line.

Applying Lemma~\ref{tran_prob_barrier}, we have
\begin{align*}
p_{\gamma}\left(0,nt,\sqrt{n}x,\sqrt{n}y \right) &= T\left(nt,(y-x)\sqrt{n}\right) - T\left(nt,(y+x)\sqrt{n}\right)\\\nonumber
 &\quad + 2\gamma \sum_{j=0}^{nt} {\gamma}^j \frac{(y+x)\sqrt{n}+j}{nt-j} T(nt-j,(y+x)\sqrt{n}+j).
\end{align*}

As $E(nt,(y\pm x)\sqrt{n})\leq M^3\epsilon^{-2}n^{-1/2}+\epsilon^{-1}n^{-1}\leq Bn^{-1/2}$, we have 
$$|1-e^{\pm C_2E(nt,(y\pm x)\sqrt{n})}|\leq Bn^{-1/2}$$
provided $n\geq n_0$ is large enough. Therefore,
\begin{equation}\label{T_1}
\begin{split}
\left| T(nt,(y\pm x)\sqrt{n})-  2  (2\pi n t)^{-1/2} e^{-(y\pm x)^2/(2t)} \right|\leq &2  (2\pi n t)^{-1/2} e^{-(y\pm x)^2/(2t)}\cdot Bn^{-1/2}\\
\leq & Bn^{-1}.   
\end{split}
\end{equation}

Fix $ \delta=1/12$. Consider the range $j\in [0,(nt)^{\frac{2}{3}-\delta}]$.
Since $\gamma=1- {\mu}/{\sqrt{n}}=e^{-\mu/\sqrt{n}\pm C\mu^2/n}$, 
$$\gamma^j=e^{-j\mu/\sqrt{n}\pm Cj\mu^2/n}=e^{-j\mu/\sqrt{n} }\exp(\pm B n^{-3\delta}).$$  

By \eqref{LLT}, we have
\begin{align*}
&T(nt-j,(y+x)\sqrt{n}+j)\times 2^{-1} ( 2\pi n t)^{1/2} e^{ [(y+x)\sqrt{n}+j]^2/(2nt)}\\
=& (1-j/(nt))^{-1/2}\exp \left( \frac{-j}{2 nt (nt-j) } [(y+x)\sqrt{n}+j]^2\pm C_2E(nt-j,(y+x)\sqrt{n}+j)    \right).
\end{align*}

We claim that, as $n\geq n_0$ large enough, the above is of the form $\exp(\pm B n^{-3\delta})$. To see the claim holds,
\begin{align*}
e^0=1\leq (1-j/(nt))^{-1/2}\leq (1-(n\epsilon)^{-1/3+\delta})^{-1/2}\leq \exp(Bn^{-1/3+\delta})\leq \exp(Bn^{-3\delta}).
\end{align*} 
\begin{align*}
0\leq \frac{ j}{2 nt (nt-j) } [(y+x)\sqrt{n}+j]^2\leq (nt)^{-4/3-\delta}[2Mn^{1/2}+(nt)^{2/3-\delta}]^2\leq Bn^{-3\delta}.
\end{align*}
\begin{align*}
 E(nt-j,(y+x)\sqrt{n}+j)=&\frac{((y+x)\sqrt{n}+j)^3}{(nt-j)^2}+\frac{1}{nt-j}\leq 4\frac{(2Mn^{1/2}+(nt)^{2/3-\delta})^2}{(nt)^2}+\frac{2}{nt}\\
\leq &	Bn^{-2/3-2\delta}\leq Bn^{-3\delta} .
\end{align*}

Hence the claim holds and we have
\begin{align*}
T(nt-j,(y+x)\sqrt{n}+j)=2( 2\pi n t)^{-1/2} e^{ -[(y+x)\sqrt{n}+j]^2/(2nt)}\exp(\pm Bn^{-3\delta}).
\end{align*}

Together with
\begin{align*}
\frac{(y+x)\sqrt{n}+j}{nt-j}=\frac{(y+x)\sqrt{n}+j}{nt}\exp(\pm Bn^{-3\delta}),
\end{align*}

we obtain that
\begin{align*}
&2\gamma\sum_{j=0}^{(nt)^{\frac{2}{3}-\delta}} {\gamma}^j \cdot\frac{(y+x)\sqrt{n}+j}{nt-j} T\left(nt-j,(y+x)\sqrt{n}+j\right)\\
=&\exp(\pm Bn^{-3\delta}) \frac{4 }{\sqrt{2\pi n t^3}} \sum_{j=0}^{(nt)^{\frac{2}{3}-\delta}}  (y+x+j/\sqrt{n})e^{-j\mu/\sqrt{n}-[y+x+j/\sqrt{n}]^2/(2t)}\cdot n^{-1/2}.
\end{align*}

As $x,y\in [0,M]$ and $t\in [\epsilon,\tau]$, we have $\int_0^\infty   (y+x+s)e^{-\mu s-(y+x+s)^2/(2t)} ds\leq  B$. Define $\textup{Err}'(n;\mu,\epsilon,M,\tau)$ to be 
\begin{align*}
\sup_{ x,y \in [0,M],t\in [\epsilon,\tau]}\left|\int_0^\infty   (y+x+s)e^{-\mu s-(y+x+s)^2/(2t)} ds-\sum_{j=0}^{(nt)^{\frac{2}{3}-\delta}}  (y+x+j/\sqrt{n})e^{-j\mu/\sqrt{n}-[y+x+j/\sqrt{n}]^2/(2t)}\cdot n^{-1/2}\right|.
\end{align*}

As the function $(y+x+s)e^{-\mu s-(y+x+s)^2/(2t)}$ decays exponentially, we have $$\lim_{n\to\infty}\textup{Err}'(n;\mu,\epsilon,M,\tau)=0.$$ 

In short,
\begin{align*}
&2\gamma\sum_{j=0}^{(nt)^{\frac{2}{3}-\delta}} {\gamma}^j \cdot\frac{(y+x)\sqrt{n}+j}{nt-j} T\left(nt-j,(y+x)\sqrt{n}+j\right)\\
=&\frac{4}{\sqrt{2\pi nt^3}} \int_0^\infty   (y+x+s)e^{-\mu s-(y+x+s)^2/(2t)} ds\pm B(n^{-3\delta-1/2}+n^{-1/2}\textup{Err}'(n;\mu,\epsilon,M,\tau)). 
\end{align*}

Next, we consider  $j\in [(nt)^{\frac{2}{3}-\delta},nt]$. Combining Lemma \ref{CLT} and Lemma \ref{pt_to_pt},
\begin{align*}
 \sum_{j\geq (nt)^{2/3-\delta}} {\gamma}^{j+1} \cdot\frac{2(y+x)\sqrt{n}+j}{nt-j} T\left(nt-j,(y+x)\sqrt{n}+j\right)= &\sum_{j\geq (nt)^{2/3-\delta}}\gamma^{j+1}\mathbb{P}^{0,x}_R(S_{nt}=y,N_{nt}(S)=j+1)\\
\leq & Be^{-n^{1/3-2\delta}/B}.
\end{align*}

Adding the above estimates, we conclude that
\begin{align*}
\left|p_{\gamma}\left(0,nt,\sqrt{n}x,\sqrt{n}y \right)-\frac{2}{\sqrt{n}}\rho_\mu (t,x,y)\right|\leq Bn^{-1/2}\big( n^{-1/2}+ n^{-3\delta }+ n^{1/2}e^{-n^{1/3-2\delta}/B}+ \textup{Err}'(n;\mu,\epsilon,M,\tau)\big).
\end{align*}

Thus the assertion follows.
\end{proof}


To prove Lemma \ref{var_transition}, we need to bound the local time for 2-D simple random walks. For $(x_1,x_2)\in\mathbb{Z}^2$, let $\mathbb{P}^{(x_1,x_2)}$ be the law of the 2-D simple random walk starting at $(S^1_0,S^2_0)=(x_1,x_2).$ For a 2-D path $(S^1,S^2)$, denote $\mathcal{N}_n$ as the number of visits to the origin before step $n-1$. In other words
\begin{align*}
\mathcal{N}_n(S_1,S_2)\coloneqq  \# \left\{ j\in [0,n-1]_\mathbb{Z} \ |\  \ (S^1_j,S^2_j)=(0,0) \right\}.
\end{align*}
The following lemma concerns the local time of 2-D random walks. The proof follows the argument in \cite[Chapter 20]{Rev}. We present the proof in the next section for the reader's convenience.

\begin{lem}\label{2D_local_time}
There exists a universal constant $C_6>0$ such that the following statement holds. For any $n\geq 2$ and $k\in\mathbb{N}_0 $, 
\begin{align*}
\mathbb{P}^{(0,0)}(\mathcal{N}_n\geq k)\leq C_6e^{ -   k / (C_6\log n)}.
\end{align*}
\end{lem}
We derive the conditional version of Lemma \ref{2D_local_time}.
\begin{lem}\label{2D_local_time_conditioned}
There exists a universal constant $C_7>0$ such that the following statement holds. For any $n\geq 2$, $k\in\mathbb{N}_0 $ and $(x_1,x_2),(y_1,y_2)\in\mathbb{Z}^2$, 
\begin{align*}
\mathbb{P}^{ (x_1,x_2)}(\mathcal{N}_n\geq k| (S^1_n,S^2_n)=(y_1,y_2))\leq  C_7e^{ - {k}/({C_7\log n})+C_7\log n }
\end{align*}
\end{lem}

\begin{proof}
We first consider the case $(x_1,x_2)=(y_1,y_2)=(0,0)$. As 
$$\mathbb{P}^{ (0,0)}(S^1_n=0,S^2_n=0)=T(n,0)^2\geq C^{-1}n^{-1},$$
we have
\begin{align*}
\mathbb{P}^{ (0,0)}(\mathcal{N}_n\geq k| (S^1_n,S^2_n)=(0,0))\leq Ce^{ -   k / (C_6\log n)+\log n}.
\end{align*}

Next, we consider general $(x_1,x_2)$ and $(y_1,y_2)$ in $\mathbb{Z}^2$. By conditioning on the first and the last time the random walk bridge touches the origin, 
\begin{align*}
 \mathbb{P}^{ (x_1,x_2)}(\mathcal{N}_n\geq k| (S^1_n,S^2_n)=(y_1,y_2)) \leq &\max_{1\leq j\leq n} \mathbb{P}^{ (0,0)}(\mathcal{N}_j\geq k-1| (S^1_j,S^2_j)=(0,0)).
\end{align*}
The change from $k$ to $k-1$ is necessary as $\mathcal{N}_n$ ignores the zero at the end. Then the assertion follows the result in the previous case.
\end{proof}

\begin{proof}[proof of Lemma \ref{var_transition}]
We adapt the notation that $B$ represents constants depending $ \tau,\epsilon$ and $K$ and $C$ represents universal constants. Without loss of generality, we assume $n\geq n_0$ with $n_0$ depending on $ \tau,\epsilon$ and $K$. The exact value of $n_0$ may increase from line to line.\\

We compute that
\begin{align*}
\E[p_{X}^2 (0,m,x,y)]&=\E\left[\sum_{S_m=y}\left( \prod_{i:S_i = 0}X_i\right)\Pro_{R}^{0,x}(S)\cdot \sum_{\tilde{S}_m=y}\left( \prod_{j:\tilde{S}_j = 0}X_i\right)\Pro_{R}^{0,x}(\tilde{S})\right]\\
&= \E\left[ \sum_{S_m=\tilde{S}_m=y}\left( \prod_ {i,j:S_i =\tilde{S}_j= 0}X_i X_j\right)\Pro_{R}^{0,x}(S)\Pro_R^{0,x}(\tilde{S})\right]\\
&= \sum_{S,\tilde{S}}\left( \prod_{i=j:S_i=\tilde{S}_i = 0}\E[X_i^2]\right)\left( \prod_{\substack{i:S_i=0\neq \tilde{S}_i \\ j:\tilde{S}_j=0\neq S_j} }\E[X_i]\E[X_j]\right)\Pro_{R}^{0,x}(S)\Pro_{R}^{0,x}(\tilde{S}).
\end{align*}

By the independence of $X$, we have
\begin{align*}
\E[p_{X}(0,m,x,y)]^2 &=\left[\sum_{S_m=y}\left( \prod_{i:S_i = 0}\E[X_i]\right)\Pro_{R}^{0,x} (S)\cdot \sum_{\tilde{S}_m=y}\left( \prod_{j:\tilde{S}_j = 0}\E[X_j]\right)\Pro_{R}^{0,x} (\tilde{S})\right]\\
&=\sum_{S_m=\tilde{S}_m=y}\left( \prod_{i=j:S_i=\tilde{S}_i = 0}(\E[X_i])^2 \right)\left( \prod_{\substack{i:S_i=0\neq \tilde{S}_i \\ j:\tilde{S}_j=0\neq S_j} }\E[X_i]\E[X_j]\right)\Pro_{R}^{0,x} (S)\Pro_{R}^{0,x} (\tilde{S}).
\end{align*}

Recall that $\E[X_i] =\gamma = 1 - {\mu}/{\sqrt{n}}$ and let $\sigma^2=\Var [X_i]$. Viewing two paths $(S ,\tilde{S})$ as a 2-D random walk, recall that $\mathcal{N}_m $ is the number of indices $i\in [0,m-1]_\mathbb{Z}$ such that $(S_i,\tilde{S}_i)=(0,0)$. We see that
\begin{align*}
\Var [p_{X}(0,m,x,y)] &= \sum_{S_m=\tilde{S}_m=y}\left(\E[X_i^2]^{\mathcal{N}_m} - \E[X_i]^{2\mathcal{N}_m}\right)\left( \prod_{\substack{i:S_i=0\neq \tilde{S}_i \\ j:\tilde{S}_j=0\neq S_j} }
\gamma \right)\Pro_{R}^{0,x}  (S)\Pro_{R}^{0,x}  (\tilde{S})\\
&= \sum_{S_m=\tilde{S}_m=y}\left((\sigma^2 + \gamma^2)^{\mathcal{N}_m} - \gamma^{2\mathcal{N}_m}\right)\left( \prod_{\substack{i:S_i=0\neq \tilde{S}_i \\ j:\tilde{S}_j=0\neq S_j} }
\gamma \right)\Pro_{R}^{0,x}  (S)\Pro_{R}^{0,x}  (\tilde{S})\\
&\coloneqq\textsc{I}_1+\textsc{I}_2.
\end{align*}
Here $\textsc{I}_1$ consists terms with $\mathcal{N}_m \leq L$ and $\textsc{I}_2$ contains terms with $\mathcal{N}_m >L$, with $L=(\log n)^{3}$. 

Suppose $\mathcal{N}_m \leq L$. For $n\geq n_0$ such that $|\mu|/n\leq 1/2$, and $Kn^{-\epsilon}/(1-\mu/\sqrt{n})^2\leq 1/2$,
\begin{align*}
 (\sigma^2 + \gamma^2)^{\mathcal{N}_m} - \gamma^{2\mathcal{N}_m} &\leq   \left((1 - {\mu}/{\sqrt{n}})^2 +Kn^{-\epsilon} \right)^{\mathcal{N}_m} - \left(1 - {\mu}/{\sqrt{n}}\right)^{2\mathcal{N}_m}\\
&= \left(1 - {\mu}/{\sqrt{n}}\right)^{2\mathcal{N}_m}\left[ \left( 1+\frac{Kn^{-\epsilon}}{(1-\mu/\sqrt{n})^2} \right)^{\mathcal{N}_m}-1\right]\\
&\leq B e^{B (\log n) ^3/\sqrt{n}}(\log n) ^3 n^{-\epsilon}\leq B(\log n) ^3 n^{-\epsilon}.
\end{align*}
Therefore we have
\begin{align*}
\textsc{I}_1 \leq & B(\log n) ^3 n^{-\epsilon} \sum_{S_m=\tilde{S}_m=y }\max\{1, \gamma\}^{N_m(S)+N_m(\tilde{S})}  \Pro^{0,x}_{R}(S)\Pro^{0,x}_R(\tilde{S}) =B(\log n) ^3 n^{-\epsilon} p_{\max\{1, \gamma\}}^2(0,m,x,y).
\end{align*}

From now on we assume $\mathcal{N}_m >L$. Let $\xi=\sigma^2+\gamma^2$. We claim that
\begin{equation}\label{I2}
\sum_{S_m=\tilde{S}_m=y,\mathcal{N}_m\geq L} \xi^{\mathcal{N}_m} \Pro^{0,x}_{R}(S)\Pro^{0,x}_R(\tilde{S})\leq  Be^{- (\log n)^2/B}  p(0,m,x,y)^2.
\end{equation}
The proof of \eqref{I2} is postponed to the end of this section. We now bound $\textsc{I}_2$ based on \eqref{I2}. Suppose $\mu\geq 0$ and hence $\gamma\leq 1$. Then from \eqref{I2},
\begin{align*}
\textsc{I}_2\leq & \sum_{S_m=\tilde{S}_m=y,\mathcal{N}_m\geq L} \xi^{\mathcal{N}_m} \Pro^{0,x}_{R}(S)\Pro^{0,x}_R(\tilde{S})\leq  Be^{- (\log n)^2/B}  p(0,m,x,y)^2.
\end{align*}

Next, we consider $\mu<0$. Let $M>0$ be a number to be determined. We further decompose $\textsc{I}_2$ into $\textsc{I}_2=\textsc{I}_{21}+\textsc{I}_{22}+\textsc{I}_{23}$. Here $\textsc{I}_{21}$ contains terms with $N_m(S),N_m(\tilde{S})\leq M\sqrt{n}$, $\textsc{I}_{22}$ contains terms with $N_m(S)>M\sqrt{n}$ and $\textsc{I}_{23}$ contains the rest.\\

If  $N_m(S),N_m(\tilde{S})\leq M\sqrt{n}$,
\begin{align*}
  \prod_{ \substack{i:S_i=0\neq \tilde{S}_i \\ j:\tilde{S}_j=0\neq S_j} } 	\gamma\leq \gamma^{2M\sqrt{n}}  \leq e^{-2\mu M}.
\end{align*}
Hence
\begin{align*}
\textsc{I}_{21} \leq & e^{-2\mu M} \sum_{\mathcal{N}_m \geq L,S_m=\tilde{S}_m=y }\left( (\sigma^2+\gamma^2)^{\mathcal{N}_m}-\gamma^{2\mathcal{N}_m} \right)\mathbb{P}_R^{0,x}(S) \mathbb{P}_R^{0,x} ( \tilde{S} )\\
\leq &Be^{-2\mu M-(\log n)^2/B}  p(0,m,x,y)^2 
\end{align*}
provided $n\geq n_0$. Here we have used the bound \eqref{I2}.\\
 
If  $N_m(S)>M\sqrt{n}$, by Cauchy-Schwarz, 
\begin{align*} 
2\textsc{I}_{22} \leq & \sum_{\mathcal{N}_m \geq L,S_m=\tilde{S}_m=y}\left( (\sigma^2+\gamma^2)^{\mathcal{N}_m}-\gamma^{2\mathcal{N}_m} \right)^2 \mathbb{P}_R^{0,x}(S) \mathbb{P}_R^{0,x} ( \tilde{S} )\\
+   & \sum_{S_m=\tilde{S}_m=y, N_m(S)> M\sqrt{n}} \gamma^{2N_m(S)+2N_m(\tilde{S})} \mathbb{P}_R^{0,x}(S) \mathbb{P}_R^{0,x} ( \tilde{S} ).
\end{align*}

From \eqref{I2}, the first term is bounded by $Be^{-(\log n)^2/B}p(0,m,x,y)^2$. The second term equals   
\begin{align*}
 &  p_{\gamma^2}(0,m,x,y)\sum_{S_m=y, N_m(S)\geq M\sqrt{n}} \gamma^{2N_m(S)}\mathbb{P}^{0,x}_R (S)\\
=& p_{\gamma^2}(0,m,x,y)p(0,m,x,y) \sum_{k\geq M\sqrt{n}} \gamma^{2k} \mathbb{P}^{0,x}_R (N_m=k|S_m=y)\\
\leq & Be^{-M	^2/B}  p^2_{\gamma^2}(0,m,x,y).   
\end{align*}

Here we have used Lemma \ref{pt_to_pt}. By symmetry, $\textsc{I}_{23}\leq \textsc{I}_{22}$ has the same upper bound. Putting the above estimates together, for $\mu<0$, $n\geq n_0$ large enough and any $M>0$,
\begin{align*}
\textsc{I}_2\leq  (Be^{-2\mu M-(\log n)^2/B}+Be^{-M^2/B}) p^2_{\gamma^2}(0,m,x,y).
\end{align*}

Choosing $M=(\log n)^2/(4B|\mu|),$ then 
\begin{align*}
\textsc{I}_2\leq   Be^{ -(\log n)^2/B}   p^2_{\gamma^2}(0,m,x,y).
\end{align*}

Thus the assertion follows.
\end{proof}
\begin{proof}[proof of \eqref{I2}]
\begin{align*}
\sum_{S_m=\tilde{S}_m=y,\mathcal{N}_m\geq L} \xi^{\mathcal{N}_m} \Pro^{0,x}_{R}(S)\Pro^{0,x}_R(\tilde{S})=\mathbb{E}^{ (x,x)}\left[ \xi^{\mathcal{N}_m} \mathbbm{1}_{ \mathcal{N}_m\geq L }  | (S_m,\tilde{S}_m)=(y,y) \right]p^2(0,m,x,y).
\end{align*}
Through summation by parts,
\begin{align*}
 &\mathbb{E}^{ (x,x)}\left[ \xi^{\mathcal{N}_m} \mathbbm{1}_{ \mathcal{N}_m\geq L }  | (S_m,\tilde{S}_m)=(y, y) \right] = \sum_{k\geq L} \xi^k \mathbb{P}^{ (x,x)}(\mathcal{N}_m= k| (S_m,\tilde{S}_m)=(y,y))\\
 = & (1-\xi^{-1})\sum_{k\geq L+1}\xi^{k}\mathbb{P}^{ (x,x)}(\mathcal{N}_m\geq k|  (S_m,\tilde{S}_m)=(y,y))+ \xi^L\mathbb{P}^{ (x,x)}(\mathcal{N}_m\geq L|  (S_m,\tilde{S}_m)=(y,y)).
\end{align*}

We require $n\geq n_0$ such that
\begin{align*}
2C^2_7 (\log(\tau n))^2 \leq &(\log n)^3,\ \xi\leq  e^{2Kn^{-\epsilon}},\ 2Kn^{-\epsilon}-1/(2C_7\log (n\tau))\leq  -1/(4C_7\log n).
\end{align*}
Here $C_7$ is the constant in Lemma \ref{2D_local_time_conditioned}. From Lemma \ref{2D_local_time_conditioned}, for any $k\geq L$,
\begin{align*}
&\mathbb{P}^{ (x,x)}(\mathcal{N}_m\geq k| (S^1_m,S^2_m)=(y,y))\leq     e^{-k/(C_7\log m)+C_7\log m}\leq e^{-k/(2C_7\log m) }.
\end{align*}

Hence
\begin{align*}
 \xi^L\mathbb{P}^{(x,x)}(\mathcal{N}_m\geq L| (S^1_m,S^2_m)=(y,y)) \leq  &\exp\left(L \left(2Kn^{-\epsilon}-1/(2C_7\log m) \right) \right)\\
\leq   \exp\left(  -L/(4C_7\log n)   \right)  \leq &e^{- (\log n)^2/B}.
\end{align*}

Similarly, 
\begin{align*}
 &(1-\xi^{-1})\sum_{k\geq L+1}\xi^{k}\mathbb{P}^{ (x,x)}(\mathcal{N}_m\geq k| (S^1_m,S^2_m)=(y,y)) \\ 
\leq  & CKn^{-\epsilon} \sum_{k\geq L+1} \exp\left(  -k/(4C_7\log n)   \right) \leq    B(\log n) n^{-\epsilon}e^{- (\log n)^2/B}\leq Be^{- (\log n)^2/B}.
\end{align*}
\end{proof}

\begin{rk}\label{third_moment}
Under the assumption $\mathbb{E}[ \left|X_i-\mathbb{E}[X_i]\right|^3 ]\leq  K n^{-\epsilon} $, we can show that $\mathbb{E}[p^3_X(0,m,x,y)]=p^3_{\gamma} (0,m,x,y)+o(1)  p^3_{\max\{1, \gamma^2\}} (0,m,x,y) $ through a similar argument because the local time of higher dimension random walks decays faster. In particular, we have $\mathbb{E}[(p_X (0,m,x,y)-p_\gamma (0,m,x,y)^\alpha)]=o(1)p^3_{\max\{1, \gamma^2\}} (0,m,x,y)$ for any $0<\alpha\leq 3$.
\end{rk}

\section{}
\begin{proof}[proof of Lemma \ref{tran_prob_barrier}]
For $m,x\geq 0 $, recall that $\mathbb{P}^{m,x}$ is the law of the symmetric simple random walk starting at $S_m=x$. For $(n,j)\in\mathbb{N}_0$, define 
\begin{align*}
q_{n,j}\coloneqq \mathbb{P}^{0,0}(S_n=j,\ S_\ell\neq j\ \textup{for all}\ \ell\in [0,n-1]_\mathbb{Z} ).
\end{align*}

As $n=0$, $[0,-1]_\mathbb{Z}$ is empty and $q_{0,j}= \mathbb{P}^{0,0}(S_0=j)=\delta_{0j}$. For $j\in\mathbb{N}_0$, define the generating function
\begin{align*}
F_j(s)&\coloneqq \sum_{n=0}^\infty q_{n,j}s^n. 
\end{align*} Note that for $j\geq 1 $, $F_j(s)=\mathbb{E}^{0,0}[s^{\tau(j)}]$ with $\tau(j)\coloneqq \inf\{n\geq 1\ |\ S_n=j\}$. By the strong Markov property, we have for $ j_1,j_2\geq 1 $, $F_{j_1} F_{j_2} =F_{j_1+j_2} .$ As $F_0=1$, the equality also holds for $ j_1,j_2\geq 0 $. In other words,
\begin{align}\label{coefq}
\sum_{k_1+k_2=n,\ k_1,k_2\geq 0} q_{k_1,j_1}q_{k_2,j_2}=q_{n,j_1+j_2}.
\end{align} 

Recall that for $(n,z)\in\mathbb{N}_0\times \mathbb{Z}$, $T(n,z)=\mathbb{P}^{0,0}(S_n=z).$ From \cite[Chapter 9]{Rev}, for any $n\geq j\geq 1$, 
\begin{equation}\label{qq}
q_{n,j}=\frac{j}{n} T(n,j).
\end{equation}

From the reflection principle, it is straightforward to derive that for any $n,j\geq 1$,
\begin{align*}
 \mathbb{P}^{0,0}(S_n=j,\ S_\ell\neq 0\ \textup{for all}\ \ell\in (0,n)_\mathbb{Z})=q_{n,j}.
\end{align*}

By conditioning on the value of $S_1$, for any $n\geq 1$,
\begin{align*}
\mathbb{P}^{0,0}(S_n=0,\ S_\ell\neq 0\ \textup{for all}\ \ell\in (0,n)_\mathbb{Z})=q_{n-1,1}.
\end{align*}

Now we start to compute $\mathbb{P}_R^{ 0,x } (N_n=0,S_{n}=y)$. By the reflection principle, for any $n\geq 1$ and $x,y\geq 0$, 
\begin{align*}
\mathbb{P}_R^{ 0,x } (N_n=0,S_{n}=y)=\left\{ \begin{array}{cc}
T(n,y-x)-T(n,y+x) & y\geq 1\\
( T(n-1,1-x)-T(n-1,1+x))/2 & y=0.
\end{array} \right. 
\end{align*}  

Assume $j\geq 1$. For any $n\geq 1$ and $x,y\geq 0$, $\mathbb{P}^{ 0,x } (N_n=j,S_{n}=y)$ equals
\begin{align*}
&\sum_{0\leq k_1<k_2<\dots k_{j}< n}\mathbb{P}^{0,x}(S_n=y, S_{k_i}=0\ \textup{for}\ i\in [1,j]_\mathbb{Z},\ S_{\ell}\neq 0\ \textup{for}\ \ell\in [0,n-1]_\mathbb{Z}\setminus\{k_1,k_2,\dots k_j\})\\
=&\sum_{0\leq k_1<k_2<\dots k_{j}< n} \mathbb{P}^{0,x}(S_{k_1}=0, S_{\ell}\neq 0\ \text{for}\ \ell\in [0,k_1-1]_\mathbb{Z})\times \mathbb{P}^{k_j,0}(S_n=y,\ \ S_{\ell}\neq 0\ \text{for}\ \ell\in (k_j,n)_\mathbb{Z})\\
&\qquad\qquad\qquad\qquad\qquad\qquad\qquad\qquad\qquad\qquad\times\prod_{i=1}^{j-1} \mathbb{P}^{k_i,0}(S_{k_{i+1}}=0,\ S_{\ell}\neq 0\ \text{for}\ \ell\in (k_i,k_{i+1})_\mathbb{Z}). 
\end{align*}
By the reflection and translation symmetry, 
\begin{align*}
\mathbb{P}^{0,x}(S_{k_1}=0, S_{\ell}\neq 0\ \text{for}\ \ell\in [0,k_1-1]_\mathbb{Z})=&q_{k_1,x},\\
\mathbb{P}^{k_i,0}(S_{k_{i+1}}=0,\ S_{\ell}\neq 0\ \text{for}\ \ell\in (k_i,k_{i+1})_\mathbb{Z})=&q_{k_i-k_{i-1}-1,1},\\
\mathbb{P}^{k_j,0}(S_n=y,\ \ S_{\ell}\neq 0\ \text{for}\ \ell\in (k_j,n)_\mathbb{Z})=&\left\{ \begin{array}{cc}
 q_{n-k_j,y} & y\neq 0,\\
 q_{n-k_j-1,1} & y=0.
\end{array} \right.
\end{align*} 
Hence for $y\neq 0$,
\begin{align*}
\mathbb{P}^{ 0,x } (N_n=j,S_{n}=y)=\sum_{0\leq k_1<k_2<\dots k_{j}< n} q_{k_1,x} q_{n-k_j,y} \prod_{i=1}^{j-1} q_{k_i-k_{i-1}-1,1}=q_{n-j+1,x+y+j-1}.
\end{align*}

Here we have used \eqref{coefq}. From \eqref{qq},
\begin{align*}
\mathbb{P}_R^{ 0,x } (N_n=j,S_{n}=y)=2\mathbb{P}^{ 0,x } (N_n=j,S_{n}=y)=\frac{2(y+x+j-1)}{n-j+1}T(n-j+1,y+x+j-1).
\end{align*}

 Similarly,
\begin{align*}
\mathbb{P}^{ 0,x } (N_n=j,S_{n}=0)=\sum_{0\leq k_1<k_2<\dots k_{j}< n} q_{k_1,x} q_{n-k_j-1,1} \prod_{i=1}^{j-1} q_{k_i-k_{i-1}-1,1}=q_{n-j ,x +j }. 
\end{align*}

Thus
\begin{align*}
\mathbb{P}_R^{ 0,x } (N_n=j,S_{n}=0)= \mathbb{P}^{ 0,x } (N_n=j,S_{n}=0)=\frac{   x+j }{n-j}T(n-j , x+j).
\end{align*}

\end{proof}

\begin{proof}[proof of Lemma \ref{CLT}]
During the proof we use $C$ to denote universal constants. Recall that for $(n,z)\in\mathbb{N}_0\times \mathbb{Z}$, $T(n,z)=\mathbb{P}^{0,0}(S_n=z).$ We first discuss the case $1/2\leq |z|/n \leq 1$. Under the assumption of the lemma we have $2^{-n}\leq T(n,z)\leq 1$. Hence
\begin{align*}
2^{-n-1}(2\pi n)^{1/2}e^{n/8} \leq  2^{-1}(2\pi n)^{1/2}e^{z^2/(2n)}  T(n,z)\leq 2^{ -1}(2\pi n)^{1/2}e^{n/2}.
\end{align*}
From the view of $E(n,z)\geq n/8$ as $|z|\geq n/2$, \eqref{CLT_2} follows for $C_1$ large enough.\\

In the rest of the proof, we assume $0\leq z<  n/2$. The case $- n/2< z\leq 0$ follows by the symmetry of $T(n,z)$. From the Stirling formula, for any $m\geq 1$,
\begin{align*}
m!\sim (2\pi m)^{1/2}m^m e^{-m}.
\end{align*}

More precisely, 
\begin{align*}
1\leq  (2\pi m)^{-1/2}m^{-m}e^m\cdot m!\leq e^{C/m} .
\end{align*}

Therefore as $n\geq 1$ we have
\begin{align}
\binom{n}{ (n+z)/2}&\leq 2^{n+1}(2\pi n)^{-1/2}\left(1-\frac{z^2}{n^2}\right)^{-n/2-1/2}\left( \frac{1-\frac{z}{n}}{1+\frac{z}{n}} \right)^{z/2}e^{C/n} \label{mid} ,\\
\binom{n}{ (n+z)/2}&\geq 2^{n+1}(2\pi n)^{-1/2}\left(1-\frac{z^2}{n^2}\right)^{-n/2-1/2}\left( \frac{1-\frac{z}{n}}{1+\frac{z}{n}} \right)^{z/2}e^{-C/(n-z)}.
\end{align}
For any $x\geq 1$,
\begin{align*}
 (1+1/x)^{-1}e&\leq  (1+1/x)^x\leq e,\  e^{-1} \leq  (1-1/x )^x\leq \ (1-1/x )^{-1}e^{-1}.
\end{align*}
We deduce $\left(1-\frac{z^2}{n^2}\right)^{z^2/2n} \leq e^{-z^2/(2n)} \left(1-\frac{z^2}{n^2}\right)^{-n/2}\leq 1$, $1\leq e^{z^2/(2n)}\left(1-\frac{z}{n}\right)^{z/2}\leq \left(1-\frac{z}{n} \right)^{-z^2/2n}$ and $1\leq e^{z^2/(2n)} \left(1+\frac{z}{n}\right)^{-z/2}\leq \left(1+\frac{z}{n} \right)^{z^2/2n}$. For any $0\leq  y \leq 1/2$, we have $1+y\leq  e^{Cy}$ and $1-y\geq e^{-Cy}$. Therefore
\begin{align*}
e^{ -Cz^4/n^3}  \leq e^{-z^2/(2n)} \left(1-\frac{z^2}{n^2}\right)^{-n/2}\leq  & 1\\
1 \leq e^{z^2/(2n)} \left(1-\frac{z}{n}\right)^{z/2},\  e^{z^2/(2n)}\left(1+\frac{z}{n}\right)^{-z/2} \leq  & e^{ Cz^3/n^2}.
\end{align*}
Similarly, $1\leq  \left(1-\frac{z^2}{n^2}\right)^{-1/2}\leq  e^{Cz^2/n^2}$. Combining the above, 
\begin{align*}
 e^{-C(z^4/n^3-1/(n-z))} \leq 2^{-(n+1)}(2\pi n)^{1/2} e^{z^2/(2n)}  \binom{n}{ (n+z)/2}&\leq  e^{ C(z^3/n^2+ z^2/n^2+1/n)} .
\end{align*}
When $0\leq z\leq n/2$, we have bounds $z^4/n^3\leq z^3/(2n^2)$, $1/(n-z)\leq 2/n$ and $z^2/n^2\leq z^3/n^2+1/n.$ Hence for $C_2$ large enough,
\begin{align*}
 e^{-C_2E(n,z) } \leq 2^{-(n+1)}(2\pi n)^{1/2} e^{z^2/(2n)}  \binom{n}{ (n+z)/2}&\leq  e^{ C_2E(n,z)}. 
\end{align*}
Thus \eqref{CLT} follows as $T(n,z)=2^{-n}\binom{n}{ (n+z)/2}$.
\end{proof}

\begin{proof}[proof of Lemma \ref{binominal_upper_bound}]
From \eqref{CLT_2},
\begin{align*}
T(n,z)\leq 2(2\pi n)^{-1/2}e^{-z^2/(2n)+C_2|z|^3/n^2+C_2/n}.
\end{align*} 

When $|z|\leq n/(4C_2)$, $C_2|z|^3/n^2\leq z^2/(4n)$ and the assertion follows by requiring $C_3\geq \max\{ 2(2\pi)^{-1/2}e^{C_2},4\}.$ Also, if $|z|=n$ then $T(n,z)=2^{-n}$ and \eqref{upper_bound} holds for $C_2$ large enough. \\

In the following, we assume $1/(4C_2) \leq z/n\leq 1-2/n$. Denote $a=z/n$. Rewriting \eqref{mid}
\begin{align*}
T(n,z)\leq 2 (2\pi n)^{-1/2}\left( {1-a^2} \right)^{-1/2}e^{-nI(a)+C/n} ,
\end{align*}

where
\begin{align*}
I(a)=\frac{1+a}{2}\ln (1+a)+\frac{1-a}{2}\ln (1-a).
\end{align*}
Since $I(a)$ is non-decreasing, $-n I(a)\leq -nI(1/(4C_2)) $ and $e^{-nI(a)}\leq e^{-n/C}$. As $a\leq 1- {2}/{n}$, $(1-a^2)^{-1/2}\leq C e^{C\log n }$. Hence
\begin{align*}
T(n,z)\leq C e^{-n/C+C\log  n +C_1/n}.
\end{align*}
Thus \eqref{upper_bound} follows as we take $C_2$ large enough.
\end{proof}

\begin{proof}[proof of Lemma \ref{2D_local_time}]
Recall that
\begin{align*}
\mathcal{N}_n(S^1,S^2)= \# \left\{ j\in [0,n-1]_\mathbb{Z} \ |\  \ (S^1_j,S^2_j)=(0,0) \right\}.
\end{align*}
Without loss of generality, we may assume $n,k\geq 3$. Define inductively $\rho_0\equiv 0$ and
$$\rho_i=\min\{j>\rho_{i-1}\ |\ (S^1_j,S^2_j)=(0,0) \}.$$
Then $$\{\mathcal{N}_n\geq k\}=\{\rho_{k-1}\leq n-1\}\subset\bigcap_{j=1}^{k-1}\{\rho_j-\rho_{j-1}\leq n-1\}.$$
As $\rho_1,\rho_2-\rho_1,\dots, \rho_{k-1}-\rho_{k-2}$ are i.i.d. 
\begin{align*}
\mathbb{P}^{ (0,0)}(\mathcal{N}_n\geq k)\leq \mathbb{P}(\rho_1\leq n-1)^{k-1}.
\end{align*}
By \cite[Lemma 20.1]{Rev}, there exists a universal constant $C$ such that $$\mathbb{P}(\rho_1\leq n-1) \leq 1- 1/(C\log(n-1))\leq e^{1/(C\log(n-1))}.$$
Hence
$$\mathbb{P}^{0,(0,0)}(\mathcal{N}_n\geq k)\leq e^{(k-1)/(C\log(n-1))}$$
and the assertion follows.
\end{proof}
\end{appendix}

\end{document}